\documentclass[twoside,11pt]{article}
\usepackage{jmlr2e}
\jmlrheading{}{}{}{1/14; Revised 5/14, 7/14}{}{2014 Vince Lyzinski, Donniell E. Fishkind, and Carey E. Priebe.}
\ShortHeadings{Seeded Graph Matching for Correlated Erd\H os-R\'enyi Graphs}{Lyzinski, Priebe, and Fishkind}
\firstpageno{1}
\usepackage{amsmath,amssymb,amsopn,amsfonts,pdfpages}
\usepackage{hyperref}
\usepackage{graphics}
\usepackage{graphicx}
\usepackage{bbm}
\usepackage{caption}
\usepackage{subcaption}

\newcommand{\R}{\mathbb{R}}

\newcommand{\tr}{\textup{trace}}

\newcommand{\p}{\mathbb{P}}

\begin{document}

\title{Seeded Graph Matching for Correlated Erd\H os-R\'enyi Graphs}

\author{\name Vince Lyzinski \email vlyzins1@jhu.edu \\
       \addr Human Language Technology Center of Excellence\\
       Johns Hopkins University\\
       Baltimore, MD, 21218, USA
       \AND
       \name Donniell E. Fishkind \email def@jhu.edu \\
       \name Carey E. Priebe \email cep@jhu.edu \\
       \addr Department of Applied Mathematics and Statistics\\
       Johns Hopkins University\\
       Baltimore, MD, 21218, USA}

\editor{}

\maketitle

\begin{abstract}
Graph matching is an important problem in  machine learning and pattern recognition.  Herein, we present theoretical and practical results on the consistency of
graph matching for estimating a latent alignment function between the vertex sets of two graphs,
as well as subsequent algorithmic implications when the latent alignment is partially observed.
In the
correlated
Erd\H{o}s-R\'enyi
graph setting, we prove that graph matching provides a strongly consistent estimate of the
 latent alignment in the presence of even modest correlation.  We then investigate
a tractable, restricted-focus version of graph matching, which
is only concerned with adjacency involving vertices in a partial observation of the latent
alignment; we prove that a logarithmic number of vertices whose alignment is known is sufficient for this
restricted-focus version of graph matching to yield a strongly consistent
estimate of the latent alignment of the remaining vertices.
We show how
Frank-Wolfe methodology for approximate graph matching, when there is
a partially observed latent alignment, inherently incorporates this restricted-focus~graph~matching. Lastly, we illustrate the relationship between seeded graph matching and restricted-focus graph matching by means of an illuminating example from human connectomics.
\end{abstract}

\section{Background and Overview}
The graph matching problem (GMP)---i.e.\@ finding the alignment between the 
vertices of two graphs which best preserves the structure of the graphs---has a rich and active place in the literature.  Graph matching has applications in a wide variety of disciplines, including machine learning \citep{nips3,pathsub,nips2}, computer vision \citep{vision3,vision2,vision}, pattern recognition \citep{pr1,pr2}, manifold and embedded graph alignment \citep{align,align2}, shape matching and object recognition \citep{nips4}, and MAP inference \citep{nips1}, to name a few.  

There are no efficient algorithms known for
solving graph matching exactly. Even the easier problem of just deciding if two
graphs are isomorphic is notoriously of unknown complexity \citep{GandJ, disease}.
Indeed, graph matching  is a special case of the NP-hard quadratic assignment
problem and, if the graphs are allowed to be directed, loopy, and weighted,
then graph matching is actually equivalent to the quadratic assignment problem.
Because of its practical applicability, there is a vast amount of literature
devoted to approximate graph matching algorithms; for an interesting survey of the literature, see e.g.\@ ``Thirty Years of Graph
Matching in Pattern Recognition" by \citet{30ygm}.

In the presence of a latent alignment function between
the vertex sets of two graphs, it is natural to ask how well graph
matching would mirror this underlying alignment.  In Section~\ref{er}
we describe the correlated Erd\H{o}s-R\'enyi random graph, which
provides us with a useful and natural setting to explore this question.
The correlated Erd\H{o}s-R\'enyi random graph consists of
two Erd\H{o}s-R\'enyi random graphs which share a common vertex set and a
common Bernoulli-trial probability parameter; for each pair of
vertices, there is a given  correlation between the two vertices' adjacency
in one graph and the two vertices' adjacency in the other graph.
In this manner, there is a natural latent alignment between the two graphs,
and we can then explore whether or not graph matching the two graphs
will consistently estimate this alignment.

If $\Phi:V(G_1)\mapsto V(G_2)$ is the latent alignment function between the vertex sets of two graphs, we define a vertex $v\in V(G_1)$ to be \emph{mismatched} by graph matching if there exists a solution $\psi$ to the graph matching problem
such that $\Phi(v)\neq \psi(v)$.  The graph matching problem provides a \emph{consistent} estimate of $\Phi$ if the number of mismatched vertices goes to zero in probability as $|V(G_1)|$ tends to infinity, and provides a \emph{strongly consistent} estimate of $\Phi$ if the number of mismatched vertices converges to zero almost surely as $|V(G_1)|$ tends to infinity.

The first of our main results is Theorem \ref{tha}, stated in
 Section \ref{er} and proven in Appendix~\ref{app}.  For correlated
 Erd\H{o}s-R\'enyi random graphs, under mild
assumptions, Theorem~\ref{tha}.i establishes that even very modest correlation is
sufficient for graph matching to yield a strongly consistent estimate of the
latent alignment; this expands and strengthens the
important results in \cite{sc}.  Theorem \ref{tha}.ii  provides a partial converse; for very weakly correlated graphs, we prove
that the expected number of permutations that align the graph more effectively (i.e. with fewer induced edge disagreements) than the latent alignment goes to infinity as the number of
vertices tends to infinity.  Unfortunately, since there is no known efficient
algorithm for graph matching, Theorem \ref{tha}.i  doesn't in-of-itself provide a
means of efficient graph alignment. However, it does suggest
that efficient approximate graph matching algorithms may be successful
in graph alignment when there is correlation between the graphs
above the threshold given in Theorem \ref{tha}.i.

Next, in Section \ref{seedgm}, we discuss the seeded graph matching
problem. This is a graph matching problem for which part of
the bijection between the two graphs' vertices is pre-specified and fixed,
and we seek to complete the bijection so as to minimize the
number of edge disagreements between the graphs; in our
correlated Erd\H{o}s-R\'enyi graph setting, the seeds are taken from
the existing latent alignment. Also in Section \ref{seedgm}, we
describe a restricted-focus version of the graph matching problem
in the context of seeding; this is a problem wherein we seek the
bijection between the two seeded graphs' vertices that minimizes
{\it only the number of seeded vertex to nonseeded vertex edge disagreements} between
the two seeded graphs. Restricting the focus of graph matching in
this particular fashion enables this restricted-focus graph matching
problem to be efficiently solved as a linear assignment problem,
in contrast to the algorithmic difficulty of (unrestricted) graph matching.

Our second main result is Theorem \ref{thaa}, which we state in
Section \ref{scerg} and prove in Appendix~\ref{app}. For correlated
Erd\H{o}s-R\'enyi graphs, under mild assumptions, Theorem \ref{thaa}.i  asserts
that a logarithmic number of seeds
is sufficient for restricted-focus graph matching to yield
a strongly consistent estimate of the latent alignment function. Theorem \ref{thaa}.ii  again provides a partial converse; for very weakly correlated graphs, we prove
that the expected number of permutations that align the unseeded vertices more effectively (i.e. with fewer induced seeded vertex to nonseeded vertex edge disagreements) than the latent alignment goes to infinity as the number of
vertices tends to infinity.
Now, what should we do if we want to perform graph alignment and there
are seeds, but the number of seeds is below this logarithmic threshold?
The remainder of this paper deals with that situation.

Back in the setting where there are no seeds,
an important class of approximate graph matching algorithms
utilize a Frank-Wolfe approach; the idea is more formally described
later in Section~\ref{sgm}. To briefly describe here, such methods
relax an integer programming formulation of graph matching to obtain
a continuous problem, then perform an iterative procedure in
which a linearization about the
current iterate is optimized, and the next iterate comes from a line
search between the current iterate and the linearization optimum. At the
conclusion of the iterative procedure, the final iterate is projected
to the nearest integer-valued point which is feasible as a graph match,
and this is taken as the approximate graph matching solution.
It turns out that the linear optimization done in each iteration can be
formulated as a linear assignment problem, which can be solved efficiently,
and this makes the Frank-Wolfe approach an appealing method
in terms of speed. The Frank-Wolfe approach can also be a very accurate
method for approximate graph matching as well; see \cite{brix,FAQ,path}
for Frank-Wolfe methodology and variants.

As done in \cite{FAP}, we describe in Section \ref{sgmalg} how this
Frank-Wolfe methodology for approximate
graph matching is naturally and seamlessly extended to the setting of seeded
graph matching so as to perform approximate seeded graph matching.
In analyzing Frank-Wolfe methodology for approximate seeded graph matching,
we observe in Section \ref{sgmcsn} that each Frank-Wolfe iteration
involves optimizing a sum of two terms.  Restricting this
optimization to just the first of these two terms turns out to
be precisely solving the aforementioned restricted-focus
graph matching problem, and restricting this optimization to just the
second of these two terms turns out to be precisely the Frank-Wolfe
methodology step if the seeds are completely ignored.

We conclude this paper with simulations and a real-data example from
human connectomics.  These simulations and experiments 
illuminate the relationship between seeded graph matching via Frank-Wolfe and restricted-focus graph matching via the Hungarian algorithm. 
We demonstrate that Frank-Wolfe methodology is often superior to
restricted-focus graph matching, an unsurprising result as the
Frank-Wolfe methodology merges
restricted-focus graph matching with seedless Frank-Wolfe methodology.
Perhaps more surprising, we also demonstrate the capacity for 
restricted-focus graph matching to outperform the full Frank-Wolfe methodology; in these cases, the noise in the unseeded adjacency can
actually degrade overall performance!

\section{Graph Matching, Random Graph Setting,
Main Results\label{sdgm}}

In this paper, all graphs will be simple graphs; in particular, edges
are undirected, there are no edges with a common vertex for both endpoints,
and there are no multiple edges between any pair of vertices. We will define $\mathcal{G}_n$ to be the
set of simple graphs on $n$ vertices.  If $G\in\mathcal{G}_n$,
we will denote the vertex set of $G$ as $V(G)$ and the edge set of $G$ via $E(G).$  For any $v,v' \in V(G)$,
if $v$ and $v'$ are adjacent in $G$ then this will be denoted $\{v,v'\}\in E(G)$,
and if $v$ and $v'$ are not adjacent in $G$ then this will be denoted
$\{v,v'\}\notin E(G)$. For any finite set $V$, the symbol $V \choose 2$ will
denote all of the ${n \choose 2}$ unordered pairs of distinct elements
from $V$.

\subsection{The Graph Matching Problem \label{gm}}

We now describe the graph matching problem.
Suppose $G_1$ and $G_2$ are graphs with the same number of vertices. Let
$\varPi$ denote the set of bijections $V(G_1) \rightarrow V(G_2)$.
For any $\psi \in \varPi $, the
{\it number of adjacency disagreements induced by $\psi$}, which will be
denoted $\Delta(\psi)$, is the number of vertex pairs
$\{v,v'\} \in {V(G_1) \choose 2}$ such that
$[\{v,v'\}\in E(G_1) \mbox{ and }\{\psi(v),\psi (v')\}\notin E(G_2)]$ or
$[\{v,v'\}\notin E(G_1) \mbox{ and }\{\psi(v),\psi (v')\}\in E(G_2)]$.
The {\it graph matching problem} is to find a bijection in $\varPi$
that minimizes the number of induced edge disagreements; we will
denote the set of solutions $\Psi:= \arg \min_{\psi \in \varPi} \Delta(\psi)$.  Equivalently stated, if $n:=|V(G_1)|=|V(G_2)|$, and if
$A,B \in \{ 0,1 \}^{n\times n}$ are the
respectively the adjacency matrices for $G_1$ and $G_2$,
then the graph matching problem is to minimize
$\| A -PBP^T \|_F$
over all $n$-by-$n$ permutation matrices $P$,
where $\| \cdot \|_F$ is the Frobenius matrix norm.

There are no efficient algorithms known for graph matching.
Even the easier problem of just deciding if $G_1$ is isomorphic to $G_2$
(i.e.~deciding if there is a bijection $V(G_1) \rightarrow V(G_2)$ which does
not induce any edge disagreements) is of unknown complexity \citep{GandJ,disease},
 and is a candidate for being in an intermediate class strictly between
P and NP-complete (if P$\ne$NP).  Also, the problem of minimizing
$\| A -PBP^T \|_F$ over all $n$-by-$n$ permutation matrices $P$, where
$A$ and $B$ are any real-valued matrices, is equivalent to the NP-hard
quadratic assignment problem. There are numerous approximate graph matching
algorithms in the literature; in Section \ref{sgm} we will discuss
Frank-Wolfe methodology.

\subsection{Correlated Erd\H{o}s-R\'enyi Random Graphs}\label{er}

Presently, we describe the correlated Erd\H{o}s-R\'enyi random graph;
this will provide a theoretical framework within which we will
prove our main theorems, Theorem~\ref{tha} and Theorem~\ref{thaa}.

The parameters are a positive integer $n$, a real
number $p$ in the interval $(0,1)$, and a real number $\varrho$
in the interval $[0,1]$; these parameters completely specify
the distribution. There is an underlying vertex set $V$
of cardinality $n$ which is common to two graphs; call these graphs
$G_1$ and $G_2$. For each $i=1,2$ and each pair of vertices
$\{v,v'\} \in {V \choose 2}$, let $\mathbbm{1}\{\{v,v'\}\in E(G_i)\}$
denote the indicator random variable for the event $\{v,v'\}\in E(G_i)$. For each $i=1,2$ and each pair of vertices
$\{v,v'\} \in {V \choose 2}$, the random variable
$\mathbbm{1}\{\{v,v'\}\in E(G_i)\}$ is Bernoulli$(p)$ distributed, and
they are all collectively independent except that, for each
pair of vertices $\{v,v'\} \in {V \choose 2}$, the
variables $\mathbbm{1}\{\{v,v'\}\in E(G_1)\}$ and $\mathbbm{1}\{\{v,v'\}\in E(G_2)\}$
have Pearson product-moment correlation coefficient $\varrho$.
At one extreme, if $\varrho$ is $1$, then $G_1$ and $G_2$ are equal, almost
surely, and at the other extreme, if $\varrho$ is $0$, then $G_1$ and
$G_2$ are independent. After $G_1$ and $G_2$ are thus realized,
their vertices are (separately) arbitrarily relabeled, so that we
don't directly observe the {\it latent alignment function (bijection)} $\Phi:V(G_1)
\rightarrow V(G_2)$ wherein, for all $v \in V(G_1)$, the vertices
$v$ and $\Phi(v)$ were corresponding vertices across the graphs
before the relabeling (i.e., the same element of $V$).

If $G_1$ is graph matched to $G_2$, to what extent will the graph match provide a consistent estimate of
the latent alignment function?
The following Theorem is our first main result.
We will be considering a sequence of random correlated
Erd\H{o}s-R\'enyi graphs with $n=1$, then $n=2$, then $n=3\ldots$, and
the parameters $p$ and $\varrho$ are each functions of $n$; i.e.~$p:=p(n)$ and
$\varrho:=\varrho(n)$.
In this paper, when we say
a sequence of events holds {\it almost always}, we mean that, with
probability~$1$, all but a finite number of the events hold.

\begin{theorem} \label{tha}
Suppose there exists a fixed real number $\xi_1<1$ such that
$p \leq \xi_1$. Then there exists fixed positive real numbers
$c_1,c_2,c_3,c_4$ (depending only on the value of $\xi_1$) such that:\\
$i)$ If $\varrho \geq c_1
\sqrt{\frac{\log n}{np}}$ and
$p \geq c_2 \frac{\log n}{n} $ then almost always $\Psi=
\{ \Phi \} $, and\\
$ii)$ If $\varrho \leq c_3 \sqrt{\frac{\log n}{n}}$ and
$p \geq c_4 \frac{\log n}{n} $ then $ \lim_{n \rightarrow \infty} \mathbb{E} |
\left \{ \psi \in \varPi : \Delta(\psi) < \Delta(\Phi)
\right  \}| =\infty$.
\end{theorem}
For proof of Theorem \ref{tha}, see Appendix \ref{app}.
\bigskip

Note that Theorem \ref{tha}.i  establishes the strong consistency of the graph matching estimate of the latent alignment function in the presence of even modest correlation between $G_1$ and $G_2$.  This theorem is a strengthening and an extension of
the pioneering work on de-anonymizing networks in
\cite{sc}, wherein the author's proved a weaker version of Theorem \ref{tha}.i  in a
sparse setting (in particular they require both $p$ and $\varrho$ to converge to $0$ at rate $p\varrho^3=O(log(n)/n)$).  Note that range of values of $p$ for which Theorem \ref{tha}.i applies includes both the sparse and the dense regimes.

Because there is no known efficient algorithm for
graph matching, Theorem \ref{tha}.i  does not directly provide a practical means of
computing the latent alignment function. But it does hold out the
hope that a good graph matching heuristic might be effective in
approximating the latent alignment function for various classes of graphs.

When proving Theorems \ref{tha} and \ref{thaa}, it will be useful for us to observe an
equivalent way to formulate correlated Erd\H{o}s-R\'enyi graphs.
For all pairs of vertices
$\{v,v'\} \in {V \choose 2}$, the indicator random
variables $\mathbbm{1}\{\{v,v'\}\in E(G_1)\}$ are independently distributed
Bernoulli$(p)$ and then (independently for the different
pairs $v,v'$), conditioning on $\mathbbm{1}\{\{v,v'\}\in E(G_1)\}=1,$ we let
$\mathbbm{1}\{\{v,v'\}\in E(G_2)\}$ be distributed Bernoulli$(p+\varrho (1-p))$ and,
conditioning on $\mathbbm{1}\{\{v,v'\}\in E(G_1)\}=0$, we let
$\mathbbm{1}\{\{v,v'\}\in E(G_2)\}$ be distributed Bernoulli$(p(1-\varrho))$.
It is an easy exercise to verify that as such, for each
$\{v,v'\} \in {V \choose 2}$, it holds that $\mathbbm{1}\{\{v,v'\}\in E(G_2)\}$ is
distributed Bernoulli$(p)$, and that the correlation of
$\mathbbm{1}\{\{v,v'\}\in E(G_1)\}$ and $\mathbbm{1}\{\{v,v'\}\in E(G_2)\}$ is $\varrho$, as desired.

\subsection{Seeded Graph Matching, Restricted-focus Graph Matching \label{seedgm}}

Continuing with the setting from Section \ref{gm}, suppose that we are
also given a subset $U_1 \subseteq V(G_1)$ of {\it seeds}
and an injective {\it seeding function}
$\phi:U_1 \rightarrow V(G_2)$, say that $U_2 \subseteq V(G_2)$
is the image~of~$\phi$. Let $\varPi_\phi$ denote the
set of bijections $\psi : V(G_1) \rightarrow V(G_2)$ such that
for all $u \in U_1$ it holds that $\psi(u)=\phi(u)$.
As before, for any bijection $\psi \in \varPi_\phi$, the
{\it number of adjacency disagreements induced by $\psi$}, which will be
denoted $\Delta(\psi)$, is the number of vertex pairs
$\{v,v'\} \in {V(G_1) \choose 2}$ such that
$[\{v,v'\}\in E(G_1) \mbox{ and }\{\psi(v),\psi(v')\}\notin E(G_2)]$ or
$[\{v,v'\}\notin E(G_1) \mbox{ and }\{\psi(v),\psi(v')\}\in E(G_2)]$.
The {\it seeded graph matching problem} is to find a bijection in
$\varPi_\phi$ that minimizes the number of induced edge disagreements;
as before, we will denote the set of solutions $\Psi:= \arg \min_{\psi \in \varPi_\phi} \Delta(\psi)$. Equivalently stated, suppose without
loss of generality that $U_1=U_2=\{v_1,v_2,\ldots,v_s\}$,
and that for all $j=1,2,\ldots,s$, $\phi(v_j)=v_j$; with $A$ and
$B$ denoting the adjacency matrices for $G_1$ and $G_2$ respectively,
the seeded graph matching problem
is to minimize $\|A - (I \oplus P)B(I \oplus P)^T\|_F$
over all $m$-by-$m$ permutation matrices $P$, where $m:=|V(G_1)|-s$, and
$\oplus$ is the direct sum, and $I$ is the $s$-by-$s$ identity matrix.

Like graph matching, there are no efficient algorithms known for
seeded graph matching; in fact, seeded graph matching is at least as
difficult as graph matching. In Section \ref{sgmalg} we discuss
how Frank-Wolfe methodology extends to provide
efficient approximate seeded graph matching.

We now present a restricted version of seeded graph matching
which  is efficiently solvable, in contrast to graph matching
and seeded graph matching.
Let $W_1:=V(G_1) \backslash U_1$ denote the nonseeds in $V(G_1)$.
For any $\psi \in \varPi_\phi$, let $\Delta_R(\psi)$ denote the number of pairs
$(w,u) \in W_1 \times U_1$ such that
$[\{w,u\}\in E(G_1) \mbox{ and }\{\psi(w),\psi(u)\}\notin E(G_2)]$ or
$[\{w,u\}\notin E(G_1) \mbox{ and }\{\psi(w),\psi(u)\}\in E(G_2)]$.
The {\it restricted-focus seeded graph matching problem} (RGM) is to find
a bijection in $\varPi_\phi$ which minimizes
such seed-nonseed adjacency disagreements; denote the set of solutions
$\Psi_R := \arg \min_{\psi \in \varPi_\phi} \Delta_R(\psi)$.
Equivalently stated, if the adjacency matrices for $G_1$ and $G_2$ are
respectively partitioned as
$$A= \begin{pmatrix} A_{11}& A_{21}^T \\ A_{21} & A_{22} \end{pmatrix}, \text{ and }
B=\begin{pmatrix} B_{11}& B_{21}^T \\ B_{21} & B_{22} \end{pmatrix} $$
where $A_{21},B_{21} \in \mathbb{R}^{|W_1| \times |U_1|}$
each represent the adjacencies between the nonseed vertices and the seed vertices
(and the seed vertices are ordered in $A_{11}$ conformally to
$B_{11}$), then finding a member of $\Psi_R$ is accomplished
by minimizing $\| A_{21} - P B_{21} \|_F$ over all
$|W_1|\times |W_1|$ permutation matrices $P$. Expanding,
\begin{align} \label{tkc}
\| A_{21} - P B_{21}\|_F^2  & =   \textup{trace} (A_{21} - P B_{21})^T
(A_{21} - P B_{21}) \notag\\ 
&  =  \textup{trace} A_{21}^TA_{21}-
\textup{trace} A_{21}^T P B_{21} - \textup{trace} B_{21}^TP^TA_{21}
+\textup{trace} B_{21}^TP^TPB_{21} \notag  \\   
&  = 
\| A_{21}\|_F^2 + \| B_{21}\|_F^2 -2 \cdot \textup{trace} \left ( P^T(A_{21} B_{21}^T)\right )  \ \   ,
\end{align}
thus finding a member of $\Psi_R$ is accomplished
by maximizing trace$P^TA_{21}B_{21}^T$ over all
$|W_1|\times |W_1|$ permutation matrices $P$. This is a linear
assignment problem and can be exactly solved in  $O(|W_1|^3)$ time with the
Hungarian Algorithm \citep{edmo,kuhn}. So, whereas finding a member
of $\Psi$ is intractable, finding a member of $\Psi_R$ can done
efficiently. An important question is how well $\Psi_R$
approximates $\Psi$.  Slightly abusing notation, we shall refer to both the restricted-focus graph matching problem and the associated algorithm for exactly solving it by RGM.

\subsection{Seeded, Correlated Erd\H{o}s-R\'enyi Graphs \label{scerg}}

Seeded, correlated Erd\H{o}s-R\'enyi graphs are
correlated Erd\H{o}s-R\'enyi graphs $G_1$ and $G_2$ where part of the
latent alignment function is observed; specifically, there is a subset
of seeds $U_1 \subseteq V(G_1)$ such that $\Phi$ is known on $U_1$.
If we take $\phi$ to be the restriction of $\Phi$ to $U_1$ and
we run RGM,
we may hope that  $\Psi_R=\{  \Phi \}$; if this hope is true
then we are provided an efficient means of computing
the latent alignment function.

The next theorem is another of our main results. We will
be considering a sequence of random correlated Erd\H{o}s-R\'enyi graphs
where the number of nonseed vertices is $m=1$, then $m=2$, then $m=3\ldots$, and the
number of seeds $s$ is a function of $m$.

\begin{theorem} \label{thaa} Suppose there exists a fixed real number $\xi_2>0$
such that $\xi_2  \leq p \leq 1-\xi_2$ and $\xi_2 \leq \varrho \leq 1-\xi_2$.
Then there exists fixed real numbers $c_5,c_6>0$ (depending only on $\xi_2$)
such that:\\
$i)$ If $s \geq c_5 \log m$ then almost always  $\Psi_R=\{ \Phi \}  $, and  \\
$ii)$ If $s \leq c_6 \log m$ then $\lim_{m \rightarrow \infty} \mathbb{E}
|\{ \psi \in \varPi_\phi : \Delta_R(\psi) < \Delta_R(\Phi) \}|=\infty$.
\end{theorem}
For proof of Theorem \ref{thaa}, see Appendix \ref{app}.

Note that Theorem \ref{thaa}.i  establishes that RGM provides a strongly consistent estimate of the latent alignment in the presence of a logarithmic number of seeds.  As noted, a member of $\Psi_R$ is efficiently computable,
and thus Theorem~\ref{thaa} (unlike Theorem~\ref{tha}) directly provides a means
to efficiently recover the latent alignment bijection $\Phi$,
if there are enough seeds.

\section{The SGM Algorithm: Extending Frank-Wolfe Methodology
for Approximate Graph Matching to Include Seeds \label{sgm}  }

In the setting with no seeds, there are numerous approximate graph matching algorithms in the literature.
One such algorithm is the FAQ algorithm of \cite{FAQ},
which is an efficient, state-of-the-art approximate graph matching algorithm based on Frank-Wolfe methodology.
The algorithm's performance is empirically shown to be state-of-the-art on many benchmark problems, and when a fixed constant number of Frank-Wolfe iterations are
performed, the running time of FAQ is O$(n^3)$, where $n$ is
the number of vertices to be matched.
Moreover, if $100 \leq | V(G_1)|$ and $G_1$ is selected with a
discrete-uniform distribution (i.e.\! all possible graphs
on $V(G_1)$ are equally likely) and $G_2$ is an isomorphic copy of $G_1$
with $V(G_2)$ being a discrete-uniform random permutation of $V(G_1)$,
then the probability that FAQ (with, say, $20$ Frank-Wolfe iterations
allowed) yields the correct isomorphism is empirically observed to be very
nearly~$1$.   We choose to focus on the FAQ algorithm here because of its amenability to seeding and because it is
the simplest algorithm utilizing the Frank-Wolfe methodology while also
achieving excellent performance on many of the QAP benchmark
problems; see \cite{FAQ}.

In Section \ref{sgmalg}, we describe the SGM algorithm from \cite{FAP}, which extends the Frank-Wolfe methodology to incorporate
utilization of seeds in approximate seeded graph matching.
In Section \ref{sgmcsn}
we point out that each Frank-Wolfe iteration in SGM involves
optimizing a sum of two terms.  Restricting this optimization
to just the first of these two terms turns out to be precisely the optimization of RGM from Section
\ref{seedgm}, and restricting this optimization to just the second of
these two terms turns out to be precisely the
corresponding optimization step of FAQ
(i.e., the seeds are completely ignored).

We conclude with simulations and real data experiments that illuminate the relationship between SGM and RGM.
SGM can be superior to RGM
matching, unsurprising in that SGM makes use of the unseeded adjacency information while RGM does not.
Perhaps more surprisingly, we also demonstrate the capacity for RGM to outperform SGM in the presence of very informative seeds; in these case the unseeded connectivity is detrimental to overall algorithmic performance!

\subsection{The Frank-Wolfe Algorithm and Frank-Wolfe Methodology \label{sgm1} }

First, a brief review of the Frank-Wolfe algorithm:
The general optimization problem that the Frank-Wolfe algorithm is applied to
is maximize $f(x)$ such that $x \in S$, where $S$ is a polyhedral set
in a Euclidean space, and the function $f:S \rightarrow \R$
is continuously differentiable. The Frank-Wolfe algorithm
is an iterative procedure. A starting point $x^{(1)} \in S$ is
chosen in some fashion, perhaps arbitrarily. For $i=1,2,3,\ldots$,
a Frank-Wolfe iteration consists of maximizing the first order
(ie linear) approximation to $f$ about $x^{(i)}$, that is maximize
$f(x^{(i)})+ \nabla f(x^{(i)})^T(x-x^{(i)})$ over $x \in S$,
call the solution $y^{(i)}$ (of course, this is equivalent to maximizing
$\nabla f(x^{(i)})^Tx$ over $x \in S$), then $x^{(i+1)}$ is
defined to be the solution to maximize $f(x)$  over $x$  on
the line segment from $x^{(i)}$ to $y^{(i)}$. Terminate the Frank-Wolfe algorithm
when the the sequence of iterates $x^{(1)},x^{(2)},\ldots$ (or their respective objective function values) stops changing much.

Of course, the seeded graph matching problem is a
combinatorial optimization problem and, as such, the Frank-Wolfe
algorithm cannot be directly applied. The term {\it Frank-Wolfe
methodology} will refer to the approach in which the
integer constraints are relaxed so that the domain is a polyhedral set and the Frank-Wolfe algorithm can be directly applied to the relaxation and, at the
termination of the Frank-Wolfe algorithm, the fractional
solution is projected to the nearest feasible integer point.
It is this projected-to point that is adopted as an
approximate solution to the original combinatorial
optimization problem. We next describe the
SGM algorithm, which applies Frank-Wolfe methodology to the Seeded Graph
Matching Problem.

\subsection{The SGM Algorithm \label{sgmalg}}

We now describe the SGM algorithm for approximate seeded graph matching.

Suppose $G_1$ and $G_2$ are graphs, say $V(G_1)=\{ v_1,v_2,\ldots,v_n \}$
and $V(G_2)=\{ v'_1,v'_2,\ldots,v'_n \}$, and let $A$ and $B$ be
the respective adjacency matrices of $G_1$ and $G_2$.
Suppose without loss of generality that
$U_1= \{ v_1,v_2,\ldots, v_s \}$ are seeds, and the
seeding function $\phi : U_1 \rightarrow V(G_2)$ is
given by $\phi (v_i)=v'_i$ for all $i=1,2,\ldots,s$.
Denote the number of nonseed vertices $m:=n-s$. Let $A$ and $B$
be partitioned
\[
 A = \left [  \begin{array}{cc} A_{11} & A_{21}^T \\
A_{21} & A_{22} \end{array} \right ]
 B = \left [  \begin{array}{cc} B_{11} & B_{21}^T \\
B_{21} & B_{22} \end{array} \right ]
\]
where $A_{11},\, B_{11} \in \{ 0,1 \}^{s \times s}$,
$A_{22},\,B_{22} \in \{ 0,1 \}^{m \times m}$, and
$A_{21},B_{21} \in \{ 0,1 \}^{m \times s}$.

As mentioned in Section \ref{seedgm}, the seeded graph matching
problem is precisely to   minimize
$\|A - (I \oplus P)B(I \oplus P)^T\|_F^2=
\| A \|_F^2 + \| B \|_F^2 -2 \cdot \textup{trace}A^T(I \oplus P) B
(I \oplus P)^T$ over all $m$-by-$m$ permutation matrices $P$.
Clearly, the seeded graph matching problem
 is equivalent to maximizing the quadratic function
$\textup{trace}A^T(I \oplus P) B
(I \oplus P)^T$ over all $m$-by-$m$ permutation matrices $P$.

Relax this maximization of $\textup{trace}A^T(I \oplus P) B
(I \oplus P)^T$ over all $m$-by-$m$ permutation matrices~$P$
to the maximization of $\textup{trace}A^T(I \oplus P) B
(I \oplus P)^T$ over all $m$-by-$m$ {\it doubly stochastic} matrices~$P$
(which form a polyhedral set), and then the Frank-Wolfe algorithm can be
applied directly to the relaxation. Simplification
yields the objective function
\begin{eqnarray}  \label{trace}f(P)  & = &
\tr A_{11}B_{11}+ \tr A_{21}^TPB_{21}+\tr A_{21} B_{21}^TP^T
+ \tr A_{22} P B_{22}P^T \\
& = &  \tr A_{11}B_{11}+ 2 \cdot \tr P^T A_{21} B_{21}^T
+ \tr A_{22}P B_{22}P^T\notag
\end{eqnarray}
which has gradient
\begin{eqnarray*}
\nabla (P) & = &  2 \cdot A_{21}B_{21}^T+2 \cdot A_{22}P B_{22} .
\end{eqnarray*}
We start the Frank-Wolfe algorithm at an arbitrarily selected doubly
stochastic $m$-by-$m$ matrix $P^{(1)}$; for convenience
we use the ``barycenter" matrix $P^{(1)}$ with all entries equal to $\frac{1}{m}$.
Then, for successive $i=1,2,\ldots$, the Frank-Wolfe iteration is to maximize
the inner product of $P$ with the gradient of $f$ at $P^{(i)}$ over all
$m$-by-$m$ doubly stochastic matrices matrices $P$; this
maximization problem is (ignoring a benign factor of $2$)
maximizing trace $P^T (A_{21}B_{21}^T + A_{22}P^{(i)} B_{22})$
over $m$-by-$m$ doubly stochastic matrices.  This
is a linear assignment problem since the optimal $P$ in this subproblem
must be a permutation matrix (by the Birkhoff-VonNeuman Theorem
which states that the $m$-by-$m$ doubly stochastic matrices are precisely
the convex hull of the $m$-by-$m$ permutation matrices), and this linear assignment
problem can be solved efficiently with the
Hungarian Algorithm in $O(m^3)$ time.
Say the optimal value of $P$ in this subproblem is
$Y^{(i)}$; then, the function $f$ on the line segment from $P^{(i)}$ to $Y^{(i)}$
is a quadratic that is easily maximized exactly, with $P^{(i+1)}$ defined as the
doubly stochastic matrix attaining this maximum.

When the Frank-Wolfe iterates
$P^{(1)},P^{(2)},P^{(3)},\ldots$ stop changing much (or
 a constant maximum of iterations are performed---we allowed 20 iterations),
then the Frank-Wolfe algorithm terminates; let the resultant approximate solution
to the relaxed problem is the doubly stochastic matrix $Q$.
The final step is to project $Q$ to the nearest $m$-by-$m$ permutation matrix.
Minimizing $\| P-Q \|_F $ over permutation matrices $P$ is again a
linear assignment problem solvable in $O(m^3)$ time;
indeed, minimizing  $$\| P-Q \|^2_F = \| P \|_F^2 -2
\mbox{trace}P^TQ + \| Q \|^2_F $$
is equivalent to
maximizing trace~$P^TQ$ over permutation matrices $P$.
This optimal permutation matrix $P$ is adopted as the approximate
solution to the seeded graph matching problem.  Specifically, the algorithm
output is the bijection $\psi : V(G_1) \rightarrow V(G_2)$ where,
for $i=1,2,\ldots,s$, \ $\psi (v_i)=v'_i$ and, for each
$i=1,2,\ldots,m$, \  $\psi (v_{s+i})=v'_{s+j}$ for the $j$ such that $P_{ij}=1$.
This Frank-Wolfe Methodology approach described above is called the {\it SGM algorithm}.

When there are no seeds, the SGM algorithm is exactly the
FAQ algorithm of \cite{FAQ};
the above development is a seamless extension of the
Frank-Wolfe methodology for approximate graph matching
when there are no seeds to Frank-Wolfe methodology for
approximate seeded graph matching.

The running time for the SGM algorithm, like for the FAQ algorithm,
is $O(n^3)$. This is because of the linear assignment problem
formulation and the use of the Hungarian algorithm in each
Frank-Wolfe iteration, and is a huge savings over using the simplex method or an interior point method for solving the linearizations
in each Frank-Wolfe iteration.
This trick has made Frank-Wolfe methodology
a very potent weapon for efficient approximate graph matching.

\subsection{Frank-Wolfe Methodology for Approximate Seeded Graph
Matching Inherently Includes RGM  \label{sgmcsn}}

In each Frank-Wolfe iteration (described in Section \ref{sgmalg}),
the linearization which is solved is
maximize (trace $P^T A_{21}B_{21}^T +  \mbox{trace}P^T
A_{22}P^{(i)} B_{22})$ over all $m$-by-$m$
permutation matrices $P$.
Observe that if this maximization were
just over the first term trace $P^T A_{21}B_{21}^T$
then it would be precisely
solving RGM from
Section \ref{seedgm}, as per Equation (\ref{tkc}) there.
Also observe that if the maximization were
just over the second term $\mbox{trace}P^T
A_{22}P^{(i)}B_{22}$, then it would be exactly the
FAQ algorithm (ignoring all of the seeds).
In this manner, the SGM algorithm can be seen
as leveraging a combination of the information gleaned from the
nonseed-seed relationships (the ``restricted-focus term") and the nonseed-nonseed relationships (the ``FAQ term").  

Although performing RGM
is much simpler than performing SGM, and although RGM almost always produces the correct graph alignment
if there are enough seeds, nonetheless SGM may perform
substantially better when there aren't enough seeds. Indeed, as
noted, SGM merges RGM with FAQ, and thus utilizes the information contained in the unseeded adjacency structure.  While FAQ alone is often unable to extract out this information (see Figure \ref{fig1} below), the RGM term can steer the FAQ term in SGM, allowing it to extract the relevant signal in the nonseed--to--nonseed adjacency structure.  

The utility of this nonseeded term depends on the amount of information captured in the seed--to--nonseed adjacency.  With less informative seeds, the SGM algorithm often significantly outperforms RGM alone, as there is important signal in the unseeded adjacency which RGM discards.  However, in the presence of well chosen seeds, the seed--to--nonseed adjacency structure may contain all the relevant signal about the unknown alignment, and the unseeded adjacency information can be a nuisance (see Figure \ref{fig:hist}).  As the RGM algorithm is exactly and efficiently solvable, this points to the centrality of both selecting and quantifying ``good'' seeds.  This is a direction of future research, as we do not address the problem of intelligent seed selection at present. 

As we will see in Figure \ref{fig1}, for weakly correlated graphs, RGM can outperform SGM.  Even with poorly-chosen seeds, the noise in the nonseed--to--nonseed adjacency structure can outweigh the relevant signal, and the performance of SGM is harmed by including this extra nuisance information.  This further highlights the utility of RGM in real data applications, where the correlation between graphs can be low.

We explore the above further in Figure \ref{fig1}.  There we compare the performance of
SGM against solving RGM
for correlated Erd\H{o}s-R\'enyi graphs with $n=300$ vertices, $p=0.5$, seeding levels ranging from $s=0$ to 275, and correlation ranging from $\varrho=0.1$ to $1$.  For each value of $\varrho$ and $s$, we ran 100 simulations and  plotted the fraction of nonseeded vertices correctly matched across the graphs, with corresponding error bars of $\pm 2$ s.e.
In all cases (except $\rho=0.1$), RGM needed more seeds to perform comparably to SGM.  Indeed, with sufficiently many seeds, all available information about the unknown alignment is captured in the seed--to--nonseed connectivity, and the (exactly solvable) RGM algorithm alone is enough to properly align the graphs.

\begin{figure}
\begin{subfigure}{.25\linewidth}
\centering
\includegraphics[trim=0cm 8cm 0cm 8cm, clip=true, width=4\textwidth]{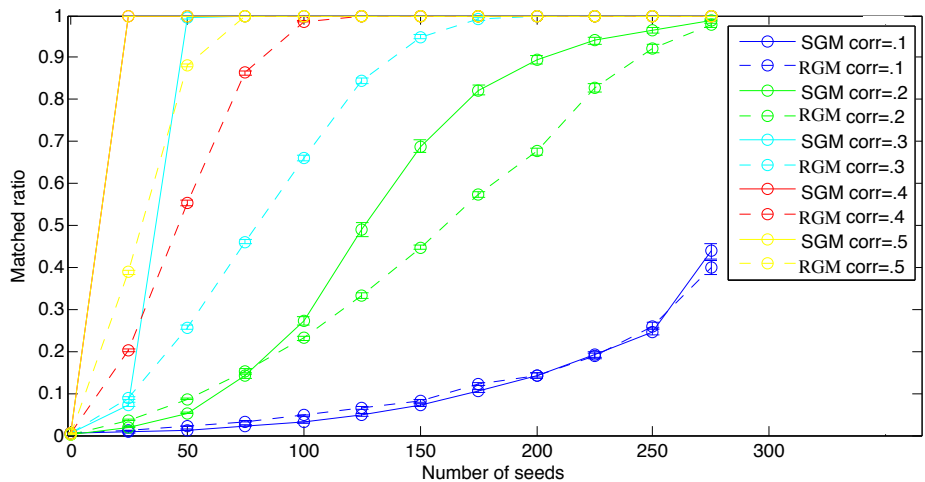}
\caption{}
\label{fig:sub6}
\end{subfigure}

\vspace{-10mm}

\begin{subfigure}{.25\linewidth}
\centering
\includegraphics[trim=0cm 8cm 0cm 8cm, clip=true, width=4\textwidth,]{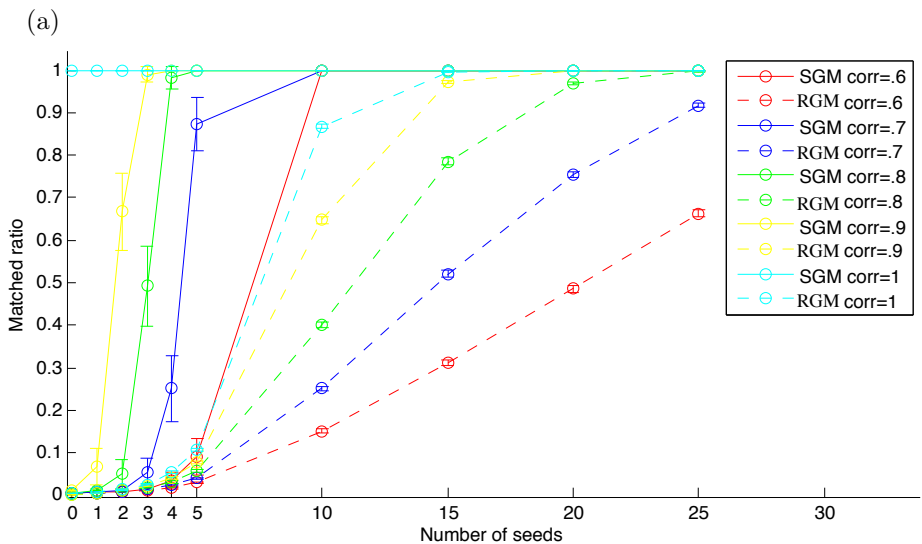}
\caption{}
\label{fig:sub7}
\end{subfigure}
\caption{Fraction of vertices correctly matched
for the SGM algorithm and for RGM, plotted
versus the number of seeds utilized, for $n=300$, $p=1/2$ and
correlation $\varrho$ varying from $0.1$ to $1$.  For each value of $\varrho$ and $s$, we ran 100 simulations and plotted the fraction of nonseeded vertices correctly matched across the graphs, with corresponding error bars of $\pm 2$ s.e.}\label{fig1}
\end{figure}
Also note the following from Figure \ref{fig1}. When there are
 no seeds,  we see FAQ (which is SGM in the absence of seeds) \
working perfectly at capturing the latent alignment function
when the two graphs are isomorphic (it bears noting that we have also observed FAQ perfectly matching when the two graphs are not isomorphic but rather *very* highly  correlated), but FAQ
does a surprisingly poor job (indeed, comparable to chance) when the correlation is even modestly less than one.
However, with seeds, SGM quickly does a very substantially better job;
indeed, the ``restricted-focus" term is steering the SGM algorithm in the proper direction!

\section{Matching Human Connectomes}

We further illuminate the relationship between SGM and RGM through a real data experiment, which will serve to highlight both the utility of RGM and the
effect of SGM's further incorporation of the unseeded adjacency information.
Our data set consists of $45$ graphs, each on $70$ vertices, these graphs constructed
respectively from diffusion tensor (DT) MRI scans of $45$ distinct healthy patients.
We have $21$ scans from the Kennedy Krieger Institute (KKI), with raw
data available at \url{http://www.nitrc.org/projects/multimodal/}, and
$24$ scans from the Nathan Kline Institute (NKI), with a description of the
raw data available at \url{ http://fcon_1000.projects.nitrc.org/indi/pro/eNKI_RS_TRT/FrontPage.html}.
All raw scans were registered to a common template and identically processed with the
MIGRAINE pipeline of \cite{wrg1}, each yielding a weighted, symmetric graph
on $70$ vertices. (All graphs can be found at
\url{http://openconnecto.me/data/public/.MR/MIGRAINE/}).
Vertices in the graphs correspond to regions in the Desikan brain atlas,
with edge weights counting the number of neural fiber bundles connecting
the regions (note that although the theory and algorithms presented earlier were for
simple graphs, they are easily modified to handle edge weights).
In addition to shedding light on the relationship between SGM and RGM,
we also explore the batch effect induced by the different medical centers
and demonstrate the capacity for seeding to potentially ameliorate this batch effect.

\begin{figure}[t]
\centering
\includegraphics[trim=0cm 3cm 0cm 2.1cm, clip=true, width=.9\textwidth]{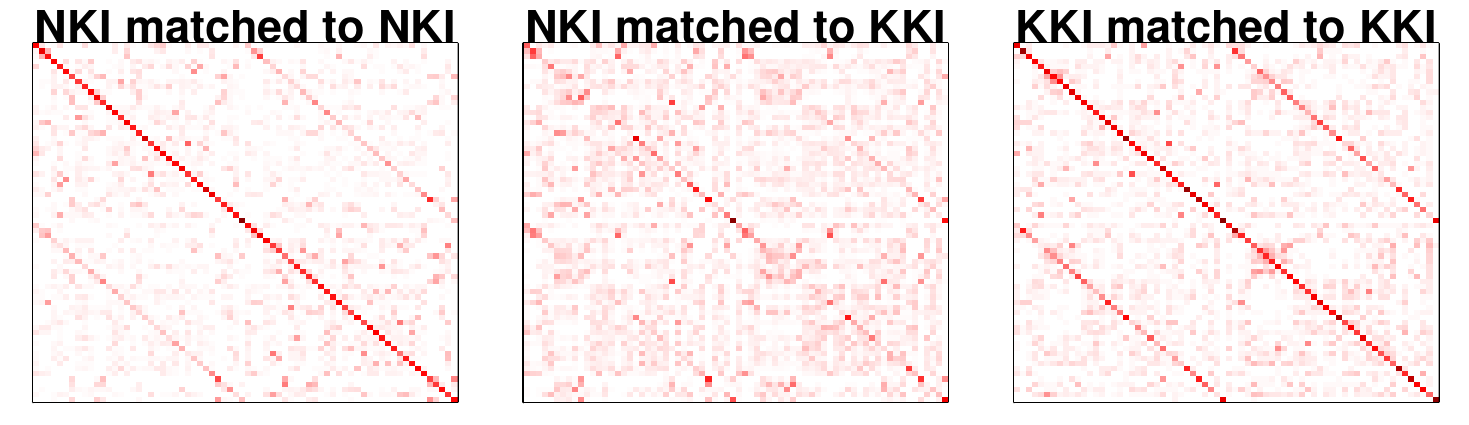}
\caption{Left: NKI to NKI matching.  Center: NKI to KKI matching. Right: KKI to KKI matching.  Each plot is a $70\times 70$ heat map with the color intensity (from white to red) representing the relative number of times vertex $i$ was match with vertex $j$ across the experiments (white denoting no matches, dark red denoting many matches).  The dark red diagonal in the left and right heat maps (as compared to the center map), indicates presence of a substantial batch effect, i.e. the correct alignment was recovered significantly better matching within medical center versus across medical center.
 Vertices 1--35 and 36--70 (as ordered) correspond to the respective brain hemispheres.}
\label{fig:batch}
\end{figure}

The pipeline which processes the scans into graphs first registers each of the graphs to a common template.
As a result, there is a canonical alignment between the vertex sets
 of these graphs (vertices corresponding to respective
 regions in the Desikan brain atlas).  How well is this alignment preserved {\it across} medical centers by the adjacency structure of the graphs alone?  Figure \ref{fig:batch} explores this question, and presents strong evidence for the existence of a batch effect (in both adjacency and geometric structure) induced by the different medical centers.  In the figure, the heatmap labeled ``KKI matched to KKI'' represents a $70\times70$ matrix, whose $i,\,j$th entry measures the
 relative number of times vertex $i$ was mapped to vertex $j$ when we ran the FAQ algorithm (i.e.\@ no seeds) over the $\binom{21}{2}$ pairs of graphs from the KKI data set.
Similarly, the ``NKI to KKI'' heat map counts the relative number of times
vertices were matched to each other when running the FAQ algorithm
over the $21\cdot24$ pairs of graphs, with one graph from each
of the KKI and NKI data sets.
The ``NKI matched to NKI'' heat map is defined similarly.
The chromatic intensity of the pixel in the $i,\,j$th entry of each heat map represents the relative frequency in which vertex $i$ was matched to vertex $j$ across the experiments, with darker red
implying more frequent and lighter red implying less frequent.
White pixels represent vertex pairs that were never matched.

Figure \ref{fig:batch} demonstrates the existence of significant signal in the adjacency structure alone (without the associated brain geometry and without seeding) for recovering the latent alignment in all three experiments.
When matching KKI to KKI, $32.8\%$ of the vertices are correctly matched on average; when matching NKI to NKI, $37.4\%$ of the vertices were correctly matched on average; when matching NKI to KKI, $9.8\%$ of the vertices were correctly matched on average (whereas chance would have matched $\approx 1.4\%$ on average).
We note that the dramatic performance difference when matching within versus across medical centers is strong evidence of the presence of a batch effect induced by the different medical centers.  Whether this batch effect is an artifact of experimental differences across medical centers (different MRI machines, different technicians, etc.) or the registration pipeline, it must be addressed before the data sets can be aggregated for use in further inference.

Also note that while much of the within medical center matching error was due to mismatching brain hemispheres (vertices 1--35 representing one
hemisphere, and vertices $36$--$70$ the other), the mismatch across medical centers appears significantly less structured.

\begin{figure}[t]
\centering
\includegraphics[trim=0cm 1.5cm 0cm 0cm, clip=true, width=1\textwidth]{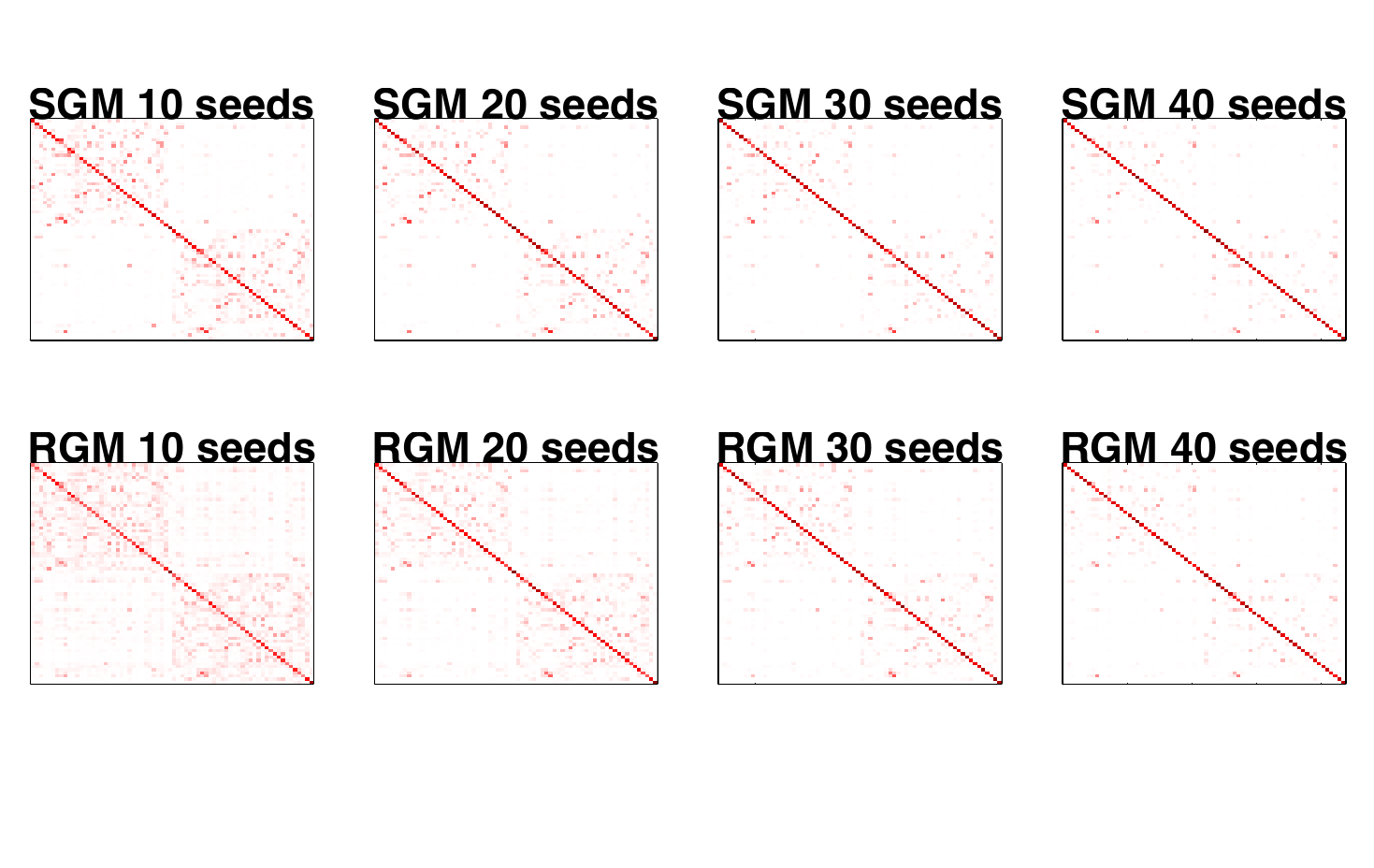}
\caption{Clockwise from top left: SGM matching the $21\cdot24$ pairs of brains, one each from the NKI and KKI data sets, using 10, 20, 30, 40 seeds; RGM matching the same set of graphs using 40, 30, 20, 10 seeds.  For each seed level, and each method we ran 100 paired MC replicates.  Each plot is a $70\times 70$ heat map with the color intensity (from white to red) representing the relative number of times vertex $i$ was matched with vertex $j$ across experiments (white denoting no matches, dark red denoting many matches).  We do not count seeded vertices as being correctly matched to each other, which would have artificially inflated the diagonal.}
\label{fig:cnssgm}
\end{figure}

Can we use seeding to ameliorate this batch effect?  In Figure \ref{fig1}, we established the capacity of seeded vertices to unearth significant signal in the adjacency structure for recovering the latent alignment function, signal which was not found without seeds.
Figure~\ref{fig:cnssgm} further demonstrates this phenomenon in our present real data setting.  We plot heat maps showing the $21\cdot24$ matchings of pairs of graphs, one each from the $NKI$ and $KKI$ data sets, for various seed levels.
For each number of seeds$=10,20,30,40,$ we ran $100$ Monte Carlo replicates (for each of SGM and RGM) for each pair of matched graphs, with each seed set chosen uniformly at random from the 70 vertices.
Clockwise, from the top left, we plot the performance of SGM with 10, 20, 30, and 40 seeds and then the performance of RGM with 40, 30, 20, and 10 seeds.
The chromatic intensity of the pixel in the $i,\,j$th entry of each heat map represents the relative frequency in which vertex $i$ was matched to vertex $j$ across the experiments
(seeded vertices are not counted as correctly matched here), with darker red implying more frequent and lighter red implying less frequent.

The figure conclusively demonstrates that seeding extracts statistically significant signal in the adjacency structure alone for correctly aligning graphs across medical center, signal that was effectively obfuscated in the absence of seeds.
While unseeded FAQ correctly matched $9.8\%$ of the vertices on average across medical centers, with $10,20,30,40$ seeds, SGM (RGM) correctly matches $49.9\%$,$68.4\%$,$78.8\%$,
$85.1\%$ ($29.9\%,53.8\%,70.7\%,80.9\%$) of the unmatched vertices on average across medical centers.  We also see that SGM outperforms RGM across all seed levels, with RGM requiring more seeds to achieve the same performance as SGM.  This is not surprising, as RGM is not utilizing any of the adjacency information amongst the unseeded vertices.

We also see
that seeding teases out additional information on the neural geometry inherent to the graphs.  For instance, with only 10 seeds, $4.3\%$ ($15.1\%$) vertices on average are mismatched across hemispheres by SGM (RGM).
In contrast, $43.8\%$ vertices on average were mismatched across hemispheres without seeds.  Interestingly, some vertex pairs are consistently mismatched across all seed levels.  For example, vertex 57 is matched by SGM to vertex 47 across medical centers $23.8\%,20.8\%,23.1\%,23.5\%$ of the time with $10,20,30,40$ seeds,  whereas, with no seeds, vertex 57 is matched to vertex 47 on average $10.9\%$ of the time when matching among the NKI data set and $10\%$ of the time when matching amongst the KKI graphs. Indeed,
these persistent artifacts are indicative of substantive differences across
(and within) data sets and demand further investigation.

\begin{figure}[h]
\centering
\includegraphics[trim=0cm .2cm 0cm 0cm, clip=true, width=1\textwidth]{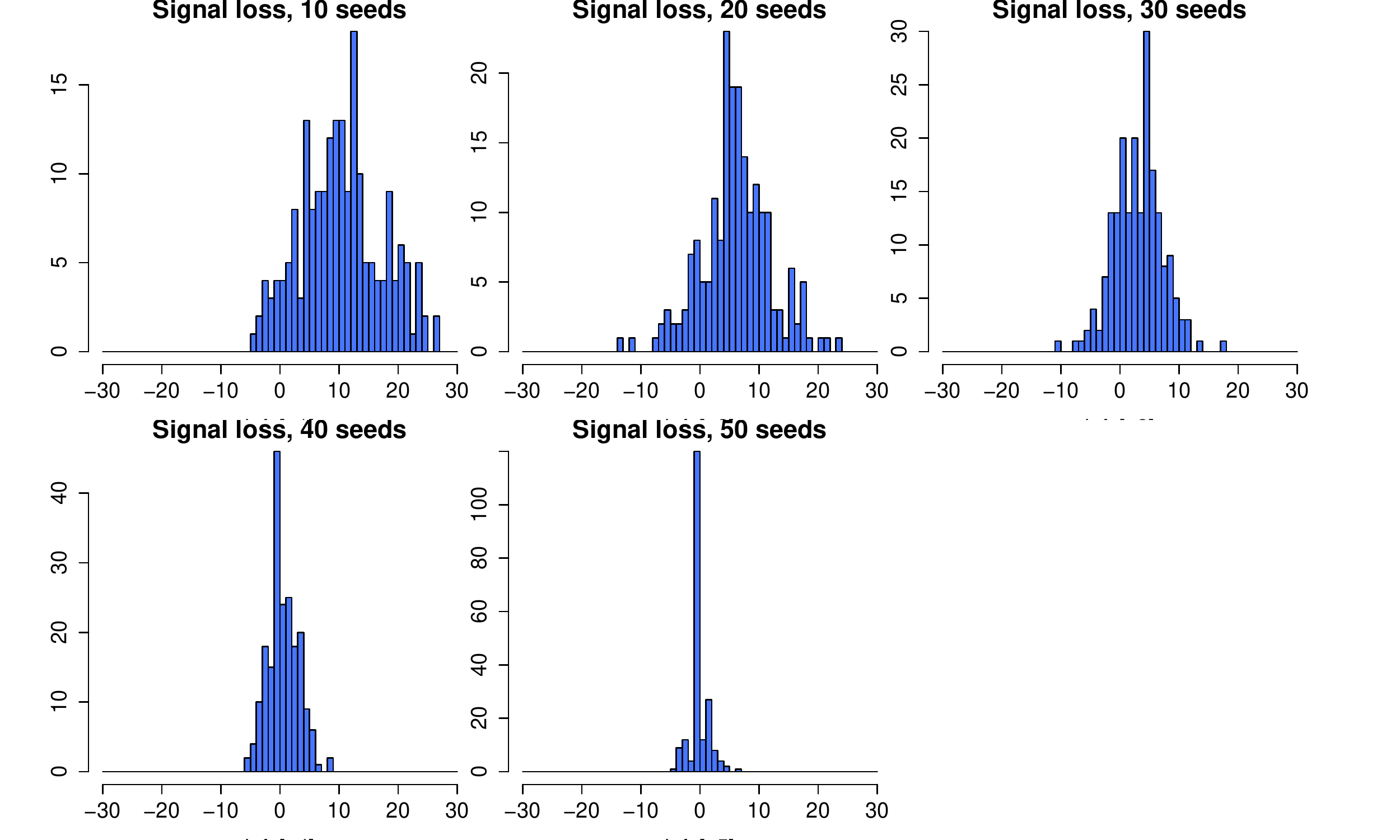}
\caption{RGM versus SGM when matching one graph each from the NKI (graph 2) and KKI (graph 7) data sets over differing seed levels.  For each of seeds$=10,20,30,40,50$ (chosen uniformly at random from the vertices), each histogram above plots 200 values of the number of vertices correctly matched by SGM minus the number of vertices correctly matched by RGM utilizing the same random seeds.}
\label{fig:hist}
\end{figure}

We have noted that, on average, SGM outperformed RGM across all seed levels.
How much of this performance gap is a function of the particular seeds chosen?  We explore this further in Figure \ref{fig:hist}.  For a pair of graphs, one each from the NKI and KKI data sets (we randomly chose graph 2 in the NKI data set and graph 7 in the KKI set---note that we see similar patterns across all tested graph pairs), we ran 200 Monte Carlo replicates of SGM and RGM seeded with the same randomly selected seeds. For each of seeds$=10,20,30,40,50$ (chosen uniformly at random from the vertices), the associated histogram plots the 200 values of the number of vertices correctly matched by SGM minus the number of vertices correctly matched by RGM.

The RGM algorithm ignores all the adjacency information amongst the unseeded vertices.
If, in Figure \ref{fig:hist}, SGM performed uniformly better than RGM at each seed level, then there is consistently relevant signal in the unseeded adjacency structure, and we should never use RGM when SGM is feasibly run.  However, we see that there are choices of seeds (at every level) for which RGM outperforms SGM.  The unseeded adjacency information is a nuisance in these cases.  As RGM is efficiently exactly solvable, this dramatically highlights the importance of intelligent seeding.  Indeed, ``good'' seeds (and hence the RGM algorithm) have the potential to capture all of the relevant adjacency structure in the graph.  While we do not pursue the question of {\it how} to select ``good" seeds here, the figure points to the centrality of this question, and we plan on pursuing active seed selection in future work.

For higher seed levels, we note that there is significantly less difference (and less variability) in the performance of SGM and RGM.
More seeds capture more information in their neighborhood structure, and the effect of the unseeded adjacency on algorithm performance is dampened.  Also, at higher seed levels the particular choice of seeds is less important, as any selection of a large number of seeds will probably contain  enough ``good" seeds to strongly align the graphs.

\section{Discussion}

Estimating the latent alignment between the vertices of two graphs is an important problem in many disciplines, and our results have both theoretical and practical implications for this problem.  Indeed, under mild assumptions, we proved the strong consistency of the graph matching problem---and its restricted focus subproblem---for estimating the latent alignment function between the vertex sets of two correlated Erd\H{o}s-R\'enyi graphs.  Although seeded graph matching is computationally hard, this result gives hope that efficient approximation algorithms will be effective in recovering the latent alignment across a broad array of graphs.

Embedded in the hard seeded graph matching problem is the tractable restricted-focus graph matching problem.  This problem is exactly solvable and also provides a strongly consistent estimator of the latent alignment.  While full seeded graph matching often out performs this restricted focus variant, we demonstrated the capacity for the restricted-focus subproblem to also outperform the full matching.  The relation between the two approaches hinges on the information contained in the seeded vertices.  If the seeds capture the adjacency structure of the graph, then the restricted-focus subproblem can benefit by {\it not} including the unseeded adjacency information, and we demonstrate this phenomenon in both real and simulation data.  This points to the primacy of intelligently seeding in graph matching, and we are working on active seeding algorithms for choosing good seeded vertices.

Even when outperformed by the full matching problem, we can still use the restricted-focus problem to extract signal in the graphs that was obfuscated without seeding.  In very large, complex problems, when it may be infeasible to run the full seeded graph matching algorithm, the restricted-focus approach could be run to provide a baseline matching between the graphs.  We are presently investigating this further, as scalability of these approaches is an increasingly important demand of modern big-data.

\vspace{3mm}

\noindent
{\bf Acknowledgments:}
This work is partially supported by a
National Security Science and Engineering Faculty
Fellowship (NSSEFF), Johns Hopkins University Human Language
Technology Center of Excellence (JHU HLT COE), and the
XDATA program of the Defense Advanced Research Projects
Agency (DARPA) administered through Air Force Research Laboratory
contract FA8750-12-2-0303. The line of research here
was suggested by Dr. \!Richard Cox, director of the JHU HLT COE.  We would like to thank Joshua Vogelstein and William Gray Roncal for providing us with the data for Section 5, and Daniel Sussman for his insightful suggestions and comments throughout.  Lastly, we would like to thank the anonymous referees whose comments greatly improved the present draft.

\bibliographystyle{plain}
\bibliography{sgmbib}

\appendix
\section{Proofs of Theorems \ref{tha} and \ref{thaa}}
\label{app}

Theorem \ref{tha} is proved in Sections \ref{tja}, \ref{thf}, and \ref{spt},
and these three subsections are a continuation one of the other.
Theorem \ref{thaa} is proved in Sections \ref{tjd}, \ref{tjf}, and
\ref{tka} and these three subsections are a continuation one of the other. Interestingly, the underlying methodology for proving Theorem~\ref{tha} is very similar (but with notable differences)
to the methodology for proving Theorem \ref{thaa}.
We begin with some results that will subsequently
be used in the proof of Theorems \ref{tha}
and \ref{thaa}.

\subsection{Supporting Results \label{supp}}

The next result, Theorem \ref{thg}, is from \cite{alon}, in the form found in
\cite{kim}.

\begin{theorem} Suppose random variable $X$ is a function of
$\eta $ independent Bernoulli$(q)$ random~variables such that changing
the value of any one of the Bernoulli random variables
changes the value of $X$ by at most $2$. For any
$t:0 \leq t<\sqrt{\eta q(1-q)}$, we have  ${\mathbb P} \left [
|X- \mathbb{E}X| > 4t \sqrt{\eta q(1-q)} \right ] \leq 2 e^{-t^2}$. \label{thg}
\end{theorem}

\noindent The next result, Theorem \ref{thi},
is a Chernoff-Hoeffding bound which is Theorem 3.2 in
\cite{fcci}.

\begin{theorem} Suppose $X$ has a Binomial$(\eta ,q)$ distribution. Then
for all $t \geq 0$ it holds that
$$\mathbb{P} \left [ X -\mathbb{E}X \geq t   \right ] \leq \text{exp}\left\{\frac{-t^2}{2
\eta q+2t/3}\right\}.$$ \label{thi}
\end{theorem}

For any $r,q \in (0,1)$, define $H(r,q):=r \log \left ( \frac{r}{q} \right ) +  (1-r) \log \left ( \frac{1-r}{1-q} \right )$. This is the Kullback-Leibler divergence
between binomial random variables with respective success probabilities $r$ and $q$.
We will later use the following rough lower bound estimate of a binomial tail probability:
\begin{proposition}
\label{P:bin}
Suppose $X$ has a Binomial($\eta ,q$) distribution, and suppose that $0<q<r<1-\frac{1}{\eta }$ for a
real number $r$. Then
$$\p(X\geq \eta r)\geq \frac{\sqrt{\pi}}{e^3}  \cdot
\sqrt{\frac{(1-r)}{r}}\eta ^{-1/2}q\cdot e^{-\eta  H(r,q)}.$$
\end{proposition}
\noindent
{\bf Proof of Proposition \ref{P:bin}}: We compute and bound
\begin{align*}
\p&(X\geq \eta r)  \geq  \p(X=\lceil \eta r\rceil)=\binom{\eta }{\lceil \eta r\rceil}q^{\lceil  \eta r\rceil}(1-q)^{\eta -\lceil \eta r\rceil}\\
&\geq  \frac{\sqrt{2\pi}}{e^2}q^{\lceil \eta r\rceil}(1-q)^{\eta -\lceil \eta r\rceil}\frac{\eta ^{\eta +0.5}}{\lceil \eta r\rceil^{\lceil \eta r\rceil+0.5}(\eta -\lceil \eta r\rceil)^{\eta -\lceil \eta r\rceil+0.5}}\\
&= \frac{\sqrt{2\pi}}{e^2}q^{\lceil \eta r\rceil}(1-q)^{\eta -\lceil \eta r\rceil}\frac{\eta ^{\eta +0.5}}{(\eta r)^{ \eta r+0.5}(\eta - \eta r)^{\eta - \eta r+0.5}}\cdot \frac{(\eta r)^{ \eta r+0.5}(\eta - \eta r)^{\eta - \eta r+0.5}}{\lceil \eta r\rceil^{\lceil \eta r\rceil+0.5}(\eta -\lceil \eta r\rceil)^{\eta -\lceil \eta r\rceil+0.5}},
\end{align*}
where the inequality in the second display line follows from Stirling's formula. Now,
\begin{eqnarray*}
 \frac{(\eta r)^{ \eta r+0.5}(\eta - \eta r)^{\eta - \eta r+0.5}}{\lceil \eta r\rceil^{\lceil \eta r\rceil+0.5}(\eta -\lceil \eta r\rceil)^{\eta -\lceil \eta r\rceil+0.5}} &=& \frac{(r)^{ \eta r+0.5}(1- r)^{\eta - \eta r+0.5}}{\left(\frac{\lceil \eta r\rceil}{\eta }\right)^{\lceil \eta r\rceil+0.5}\left(1-\frac{\lceil \eta r\rceil}{\eta }\right)^{\eta -\lceil \eta r\rceil+0.5}}\\
&\geq & \left(\frac{1}{\frac{\lceil \eta r\rceil}{ \eta r}}\right)^{\eta r+0.5}(1-r)^{\lceil \eta r\rceil-\eta r}\\
&\geq & \left(\frac{1}{1+\frac{1}{ \eta r}}\right)^{\eta r+0.5}(1-r)
\ \ \ \ \geq \ \ \ \ \frac{1}{e\sqrt{2}} (1-r).
\end{eqnarray*}
Combining the above, we obtain
\begin{eqnarray*}
\p(X\geq \eta r) &\geq & \frac{\sqrt{\pi}}{e^3}  \cdot \frac{1-r}{r^{1/2}(1-r)^{1/2}} \eta ^{-1/2}q^{\eta r+1}(1-q)^{\eta -\eta r}\frac{\eta ^{\eta }}{(\eta r)^{ \eta r}(\eta - \eta r)^{\eta - \eta r}}\\
&= & \frac{\sqrt{\pi}}{e^3}  \cdot \sqrt{\frac{1-r}{r}}\eta ^{-1/2}q\cdot e^{-\eta H(r,q)},
\end{eqnarray*}
as desired. $\blacksquare$

\subsection{Overall Argument of the Proof for Theorem \ref{tha}, Part i \label{tja}}

It is notationally convenient to assume without loss of generality
that the correlated Erd\H{o}s-R\'enyi graphs $G_1$ and $G_2$ are on
the same set of $n$ vertices $V$ and we do {\bf not}
relabel the vertices. Let $\varPi$ denote the set of bijections
$V \rightarrow V$; here, the identity function $e \in \varPi$ is
the latent alignment bijection~$\Phi$.
For any $\psi \in \varPi$, 
\begin{align*}\Delta^+(G_1,G_2,\psi)&:=| \left \{  \{ v,v' \} \in {V \choose 2}
\mbox{ such that }
\{v,v'\}\notin E(G_1) \mbox{ and }  \{\psi(v),\psi(v')\}\in E(G_2) \right \}| ,\\
\Delta^-(G_1,G_2,\psi)&:=| \left \{  \{ v,v' \} \in {V \choose 2}
\mbox{ such that }
\{v,v'\}\in E(G_1)  \mbox{ and } \{\psi(v),\psi(v')\}\notin E(G_2)  \right \}| ,\\
\Delta^{0+}(G_1,G_2,\psi)&:=| \bigg \{ \{ v,v' \} \in {V \choose 2}
\mbox{ such that }
\{v,v'\}\notin E(G_1),\,\{\psi (v),\psi (v')\}\in E(G_1),\, \\
&\hspace{25mm}\{\psi (v),\psi (v')\}\notin E(G_2) \bigg \} |
, \\
 \Delta^{0-}(G_1,G_2,\psi)&:=| \bigg \{ \{ v,v' \} \in {V \choose 2}
\mbox{ such that }
\{v,v'\}\in E(G_1),\,\{\psi (v),\psi (v')\}\notin E(G_1),\, \\
&\hspace{25mm}\{\psi (v),\psi (v')\}\in E(G_2) \bigg \} |,\\
\Delta (G_1,G_2,\psi)&:=\Delta^+ (G_1,G_2,\psi) + \Delta^- (G_1,G_2,\psi)
.\end{align*}

First, note that
\begin{eqnarray} \label{the}
\Delta^+ (G_1,G_1,\psi)=\Delta^- (G_1,G_1,\psi)= \frac{1}{2}
\Delta (G_1,G_1,\psi) \ \ ;
\end{eqnarray}
this is because the number of edges in $G_1$ isn't changed when
its vertices are permuted by $\psi$.

Next, note that
\begin{align} \label{thb}
\Delta(G_1,G_2,\psi)-\Delta(G_1,G_2,e)=\Delta(G_1,G_1,\psi)-2 \cdot \Delta^{0+}(G_1,G_2,\psi) -2 \cdot \Delta^{0-}(G_1,G_2,\psi) \ \ ;
\end{align}
this is easily verified by replacing
``$G_2$" in (\ref{thb}) by ``$G$", and
observing the truth of (\ref{thb}) as $G$, starting out with $G=G_1$,
is changed one edge-flip at a time until $G=G_2$.

Now, consider the event, which we shall call $\Upsilon$, that
for all $\psi \in \varPi\backslash\{e\}$,
\begin{eqnarray} \label{thc}
\Delta^{0+}(G_1,G_2,\psi) & < & \Delta^{+}(G_1,G_1,\psi) \cdot \left ( (1-p)(1-\varrho)+\frac{\varrho}{2} \right ) \mbox{  and also}\\ \label{thd}
\Delta^{0-}(G_1,G_2,\psi) & < & \Delta^{-}(G_1,G_1,\psi) \cdot \left ( p(1-\varrho)+\frac{\varrho}{2} \right ).
\end{eqnarray}
We will next show in Section \ref{thf} that, under the hypotheses of the first
 part of Theorem \ref{tha}, $\Upsilon$ almost
always happens (in other words, with probability $1$, $\Upsilon$ happens
for all but a finite numbers of $n$'s).
Then, adding (\ref{thc}) to (\ref{thd}) and using
(\ref{the}), we then obtain that almost always
$\Delta^{0+}(G_1,G_2,\psi) + \Delta^{0-}(G_1,G_2,\psi) < \frac{1}{2} \cdot
\Delta (G_1,G_1,\psi) $ for all $\psi \in \varPi\backslash\{e\}$.
Substituting this into (\ref{thb}) yields
that almost always $\Delta(G_1,G_2,\psi)>\Delta(G_1,G_2,e)$
for all $\psi \in \varPi\backslash\{e\}$, and the first part
of Theorem \ref{tha} is then proven.

\subsection{Under Hypotheses of Theorem \ref{tha}, Part i, $\Upsilon$ Occurs Almost Always\label{thf}}

For any $k \in \{ 1,2,\ldots, n \}$, let  $\varPi(k)$ denote the set of bijections
in $\varPi$ such that the number of non-fixed-points of the bijection is exactly $k$; that is,
$\varPi(k):= \{ \psi \in \varPi:
| \{ v \in V: \psi(v) \ne v \}| = k \}$. A simple upper bound for $|\varPi(k)|$ is
$|\varPi(k)| \leq {n \choose k} k! = n(n-1)(n-2) \cdots (n-k+1)  \leq n^k$.

Just for now, let $k \in \{ 1,2,\ldots, n \}$ be chosen, and let $\psi
\in \varPi(k)$ be chosen.
Denoting $T(\psi):=\left\{ \{ v,v' \} \in {V \choose 2} \mbox{ such that } v=\psi(v'), \, v'=\psi(v)\}\right\}$, we have that
the random variable $\Delta(G_1,G_1,\psi)$ is a function
of the $\eta := {k \choose 2} + (n-k)k-|T(\psi)|$ independent Bernoulli$(p)$
random variables 
$$ \{ \mathbbm{1}\{\{v,v'\}\in E(G_1)\} \}_{ \{ v,v' \} \in {V \choose 2}\backslash T(\psi)\ : \
\psi(v)\ne v \textup{ or }\psi(v') \ne v'  }$$ and note that the hypotheses
of Theorem \ref{thg} are satisfied, hence for the choice of
$t=\frac{1}{20} \sqrt{\eta p(1-p)}$ in Theorem \ref{thg} we obtain that
\begin{eqnarray} \label{tia}
\mathbb{P} \left [ | \Delta(G_1,G_1,\psi)- \mathbb{E}\Delta(G_1,G_1,\psi) |>
\frac{1}{5} \eta p(1-p) \right ] \leq 2e^{- \eta  p (1-p)/400 }.
\end{eqnarray}
Also note that
\begin{eqnarray*}
\Delta(G_1,G_1,\psi)=\sum_{\substack{\{ v,v' \} \in {V \choose 2}\backslash T(\psi) \\ \text{s.t.} \
\psi(v)\ne v \textup{ or }\psi(v') \ne v'}}
 \mathbbm{1}\bigg\{ \mathbbm{1}\{\{ v,v' \}\in E(G_1)\}
\ne \mathbbm{1}\{ \{ \psi(v),\psi(v') \}\in E(G_1)\}  \bigg\}
\end{eqnarray*}
is the sum of $\eta $ Bernoulli$(2p(1-p))$ random variables hence
\begin{eqnarray} \label{tib}
\mathbb{E} \Delta(G_1,G_1,\psi)=2 \eta p(1-p) .
\end{eqnarray}
Because $|T(\psi)|\leq \frac{k}{2}$, we have by
elementary algebra that $\frac{(n-2)k}{2} \leq
\eta \leq nk$. Thus, by (\ref{tia}) and (\ref{tib})
we obtain that (for large enough $n$; in the following our constants
are very conservatively chosen)
\begin{eqnarray} \mathbb{P} \left ( \frac{\Delta(G_1,G_1,\psi)}{nkp(1-p)} \not \in
[1/2,\  5/2] \right ) \leq 2 e^{\frac{-1}{1000}nkp(1-p)}
\leq 2 e^{\frac{-(1-\xi_1)}{1000}nkp}  .  \label{tie}
\end{eqnarray}

Conditioning on $G_1$, random variable $\Delta^{0+}(G_1,G_2,\psi)$ has a
$$\text{Binomial}\left ( \Delta^+(G_1,G_1,\psi),(1-p)(1-\varrho) \right )$$ distribution, and
random variable $\Delta^{0-}(G_1,G_2,\psi)$ has a $$\text{Binomial}\left (
\Delta^-(G_1,G_1,\psi),p(1-\varrho) \right )$$distribution.
Conditioning also on the event that $\frac{\Delta(G_1,G_1,\psi)}{nkp(1-p)} \in
[1/2,\  5/2]$, we apply Theorem~\ref{thi} with the value $t=\frac{\varrho}{2}
\cdot \Delta^+(G_1,G_1,\psi)$, and we use (\ref{the}) to show
\begin{eqnarray} \label{tic}
\mathbb{P} \left [  \Delta^{0+}(G_1,G_2,\psi) \geq
\Delta^+(G_1,G_1,\psi) \cdot \left (  (1-p)(1-\varrho)+ \frac{\varrho}{2} \right )
 \right ]  \leq e^{\frac{-(1-\xi_1)}{40}nkp\varrho^2},\\ \label{tid}
\mathbb{P} \left [  \Delta^{0-}(G_1,G_2,\psi) \geq
\Delta^-(G_1,G_1,\psi) \cdot \left (  p(1-\varrho)+ \frac{\varrho}{2} \right )
 \right ] \leq e^{\frac{-(1-\xi_1)}{40}nkp\varrho^2}.
\end{eqnarray}

Finally, applying (\ref{tie}), (\ref{tic}) and (\ref{tid}), the probability of $\Upsilon^C$ can be bounded using subadditivity:
\begin{eqnarray*}
\mathbb{P}(\Upsilon^C) & \leq & \sum_{k=1}^n \sum_{\psi \in \varPi(k)}
\left (  2 e^{\frac{-(1-\xi_1)}{1000}nkp} + e^{\frac{-(1-\xi_1)}{40}nkp\varrho^2}
  +   e^{\frac{-(1-\xi_1)}{40}nkp\varrho^2} \right ) \\
 & \leq &  \sum_{k=1}^n n^k \left (   2 e^{\frac{-(1-\xi_1)}{1000}nkp} +2 e^{\frac{-(1-\xi_1)}{40}nkp\varrho^2} \right )\\
 & \leq & \sum_{k=1}^n  \left (   2 e^{\frac{-(1-\xi_1)}{1000}nkp+k \log n} +2e^{\frac{-(1-\xi_1)}{40}nkp\varrho^2 + k \log n}   \right )
 \leq n \cdot \frac{4}{n^3},
\end{eqnarray*}
the last inequality holding if
$p \geq c_2 \frac{\log n}{n} $ and
$\varrho \geq c_1 \sqrt{\frac{\log n}{np}}$ for sufficiently large, for fixed
constants $c_1,c_2$. Because $\sum_{n=1}^\infty \frac{4}{n^2}<\infty$, we have by
the Borel-Cantelli Lemma that $\Upsilon$ almost always happens.
As mentioned in Section \ref{tja}, this completes
the proof of the first part of Theorem \ref{tha}. $\blacksquare$
\begin{remark}
\emph{Note that we could tighten the constants $c_1$ and $c_2$ appearing above.  Here we choose not to, instead focusing on the orders of magnitude of $\varrho$, and do not pursue exact constants further.
}\end{remark}

\subsection{Proof of Theorem \ref{tha}, Part ii\label{spt}}

We now prove the second part of Theorem \ref{tha}.

Just for now, let $\psi \in \varPi(n)$ be chosen (i.e., $\psi$ is a derangement),
and condition on $\Delta^+(G_1,G_1,\psi)=\Delta$,
where $ \frac{1}{4}n^2p(1-p) \leq \Delta \leq
\frac{5}{4}n^2p(1-p)$. The random variables
$\Delta^{0+}(G_1,G_2,\psi)$ and $\Delta^{0-}(G_1,G_2,\psi)$
are independent, and have distributions Binomial$(\Delta,q_1)$ and Binomial$(\Delta,q_2)$, respectively,
where $q_1:=(1-p)(1-\varrho)$ and $q_2:=p(1-\varrho)$.

Denoting $r_1:=q_1+\frac{\varrho}{2}$ and $r_2:=q_2+\frac{\varrho}{2}$,
and observing that, under the hypotheses of Theorem \ref{thaa}, part ii,
it holds that $r_1 < 1-\frac{1}{\Delta}$ and $r_2 < 1-\frac{1}{\Delta}$
we thus have by Proposition~\ref{P:bin} that (as $\frac{\pi}{e^6}>\frac{1}{200}$)
\begin{align*}
\mathbb{P} \Bigg ( \Delta^{0+}(G_1,G_2,\psi) &\geq \Delta \cdot r_1 \ \mbox{ and } \
\Delta^{0-}(G_1,G_2,\psi) \geq \Delta \cdot r_2 \Bigg )\\
& \geq
\frac{q_1q_2}{200 \Delta}\sqrt{\frac{(1-r_1)(1-r_2)}{r_1r_2}}
e^{-\Delta \cdot H(r_1,q_1) - \Delta \cdot H(r_2,q_2)}
\end{align*}

Note that we can change the inequalities ``$\geq$" in the expression
 $\mathbb{P}(\ \ \ )$ above into
 strict inequalities ``$>$"  with a harmless tweak.
An elementary calculus argument yields that
$H(x+y,y) \leq x^2/(y-y^2)$ for all
$0<y<1$ and $x \geq 0$ such that $y+2x<1$. Indeed, fixing any value for $y$,
the function value and the derivative of $H(x+y,y)$ with respect to $x$ are
both~$0$ at $x=0$, the function value and the derivative of $x^2/(y-y^2)$
with respect to $x$ are both $0$ at $x=0$, and the second derivative of
$H(x+y,y)$ is less than the second derivative of $x^2/(y-y^2)$
for all relevant $x$.  This, together with the fact
that $1-r_1=r_2$, $1-r_2=r_1$ and assuming that $\varrho $ is bounded away from $1$
(which, indeed, will turn out to be assumed),
we have that there exists a real number $c>0$ such that
\begin{align} \label{tjb}
\mathbb{P} \Bigg ( \Delta^{0+}(G_1,G_2,\psi) > \Delta \cdot &r_1 \mbox{ and }
\Delta^{0-}(G_1,G_2,\psi) > \Delta \cdot r_2 \Bigg )  \notag  \\
 &\geq  \frac{q_1q_2}{200\Delta}e^{-\varrho^2 \Delta \left (  \frac{1}{4q_1(1-q_1)} +
\frac{1}{4q_2(1-q_2)}  \right )} \notag\\
& \geq  \frac{c}{n^2} \cdot e^{-\varrho^2 n^2p \cdot \left (
\frac{1}{c \cdot p} \right ) } \notag \\
 &= \frac{c}{n^2} \cdot e ^{\frac{-\varrho^2 n^2 }{c} } \ \ . \ \ \  \ \ \ \ \  \ \
\end{align}

From (\ref{the}) and (\ref{tie}) we have that
there exists a fixed constant $c_4$ such that
if $p \geq c_4 \frac{\log n}{n}$ then,
with probability $>\frac{1}{2}$ (for $n$ large enough)
it holds that $ \frac{1}{4}n^2p(1-p) \leq \Delta^+(G_1,G_1,\psi) \leq \frac{5}{4}n^2p(1-p)$.
Thus, by (\ref{tjb}), noting again that
 $r_1+r_2=1$ and that $\Delta^+(G_1,G_1,\psi)=\frac{1}{2}\Delta(G_1,G_1,\psi)$, we have unconditionally
\begin{eqnarray}
\mathbb{P} \Bigg (
\Delta^{0+}(G_1,G_2,\psi) + \Delta^{0-}(G_1,G_2,\psi) > \frac{1}{2} \cdot
\Delta (G_1,G_1,\psi) \Bigg ) \geq \frac{c}{2n^2} \cdot e ^{\frac{-\varrho^2 n^2 }{c} } \label{tjc}
\end{eqnarray}

Next, the number of derangements $|\varPi(n)|$ satisfies $\lim _{n \rightarrow \infty}
\frac{|\varPi(n)|}{n!}=\frac{1}{e}$, thus with Stirling's formula we have that for $n$ large enough it will hold that $|\varPi(n)| \geq \left (  \frac{n}{e}   \right )^n$.
Thus, for $n$ large enough, by  (\ref{thb}) and (\ref{tjc}),
\begin{eqnarray*}
\mathbb{E} |
\left \{ \psi \in \varPi : \Delta(G_1,G_2,\psi) < \Delta(G_1,G_2,e)
\right  \}| &=& \sum_{\psi \in \varPi} \mathbb{P}
\Bigg (\Delta(G_1,G_2,\psi) < \Delta(G_1,G_2,e)  \Bigg )\\
& \geq & \sum_{\psi \in \varPi(n)} \mathbb{P}
\Bigg (\Delta(G_1,G_2,\psi) < \Delta(G_1,G_2,e)  \Bigg )\\
& \geq & \left (  \frac{n}{e}   \right )^n  \frac{c}{2n^2} \cdot e ^{\frac{-\varrho^2 n^2 }{c} } \\
&=&
\frac{c}{2n^2} \cdot e ^{\frac{-\varrho^2 n^2 }{c} +n \log n - n},
\end{eqnarray*}
so that there exists a fixed real number $c_3>0$ such that
if $\varrho \leq c_3 \sqrt{\frac{\log n}{n}}$
then it holds that $\mathbb{E} |
\left \{ \psi \in \varPi(n) : \Delta(G_1,G_2,\psi) < \Delta(G_1,G_2,e)
\right  \}| \rightarrow \infty $ as $n \rightarrow \infty$,
and the second part of Theorem~\ref{tha} is proven. $\blacksquare$
\begin{remark}
\emph{Note that we could tighten the constants $c_3$ and $c_4$ appearing above.  Here we choose not to, instead focusing on the orders of magnitude of $\varrho$, and do not pursue exact constants further.
}
\end{remark}
\subsection{Overall Argument of the Proof for Theorem \ref{thaa}, part i \label{tjd}}

The proof of Theorem \ref{thaa} is very similar in structure to the
proof of Theorem \ref{tha}. For simplicity of notation,
suppose without loss of generality that the correlated Erd\H{o}s-R\'enyi
graphs $G_1$ and $G_2$ are on the same set of $n$ vertices $V$, and we
do {\bf not} relabel the vertices. Let $\varPi$ denote the set of bijections
$V \rightarrow V$; here the identity function $e \in \varPi$ is
the latent alignment bijection. Further suppose that $V$ is partitioned into
$s$ seed vertices $U$, and $m$ nonseed vertices $W$. Let $\phi: U \rightarrow U$
be the identity function, and let $\varPi_\phi :=
\{ \psi \in \varPi : \forall u \in U \ \psi(u)=u \}$. For any $\psi \in \varPi_\phi$, define 
\begin{align*}
\Delta_R^+(G_1,G_2,\psi)&:= | \{ (w,u) \in W \times U :\{w,u\}\notin E(G_1)
\mbox{ and } \{\psi(w),u\}\in E(G_2) \} | ,\\
\Delta_R^-(G_1,G_2,\psi)&:= | \{ (w,u) \in W \times U :\{w,u\}\in E(G_1)
\mbox{ and } \{\psi(w),u\}\notin E(G_2) \} |,\\
\Delta_R^{0+}(G_1,G_2,\psi)&:= \big| \{ (w,u) \in W \times U :\{w,u\}\notin E(G_1),\, 
\{\psi(w),u\}\in E(G_1),\,\\
&\hspace{75mm}\{\psi(w),u\}\notin E(G_2)\} \big| ,\\
\Delta_R^{0-}(G_1,G_2,\psi)&:= \big| \{ (w,u) \in W \times U :\{w,u\}\in E(G_1),\, 
\{\psi(w),u\}\notin E(G_1),\,\\
&\hspace{75mm}\{\psi(w),u\}\in E(G_2)\} \big| ,\\
\Delta_R(G_1,G_2,\psi)&:=\Delta_R^+(G_1,G_2,\psi)+\Delta_R^-(G_1,G_2,\psi).\end{align*}

First note that
\begin{eqnarray} \label{tjk}
\Delta_R^+(G_1,G_1,\psi)=\Delta_R^-(G_1,G_1,\psi)= \frac{1}{2}
\Delta_R(G_1,G_1,\psi) \ \ ;
\end{eqnarray}
this can be easily verified by considering, for each $u \in U$
and for each cycle $C$ of the permutation $\psi$, the changes
of status in adjacency-to-$u$ of the vertices as the vertices of $C$ are
considered in their cyclic order. (Specifically, the number of changes along $C$
from adjacency-to-$u$ to nonadjacency-to-$u$ are equal to the number
of changes along $C$ from nonadjacency-to-$u$ to adjacency-to-$u$.)

Next, note that
\begin{eqnarray} \label{tje} \label{tjm}
 \Delta_R (G_1,G_2,\psi)-  \Delta_R (G_1,G_2,e)   = \Delta_R (G_1,G_1,\psi)
-2  \cdot \Delta_R ^{0+}(G_1,G_2,\psi)-2 \cdot \Delta_R^{0-}(G_1,G_2,\psi) \ \ ; \ \ \ \ \ \
\end{eqnarray}
this is easily verified by replacing ``$G_2$" in (\ref{tje}) with ``$G$", and observing the
truth of (\ref{tje}) as $G$, starting out with $G=G_1$, is changed
one edge-flip at a time until $G=G_2$.

Now, consider the event $\Upsilon_R $
defined as it holding that, for all $\psi \in \varPi_\phi$ besides $e$,
\begin{eqnarray} \label{tji}
\Delta_R^{0+}(G_1,G_2,\psi) & < & \Delta_R^{+}(G_1,G_1,\psi) \cdot \left ( (1-p)(1-\varrho)+\frac{\varrho}{2} \right ) \mbox{  and also}\\ \label{tjj}
\Delta_R^{0-}(G_1,G_2,\psi) & < & \Delta_R^{-}(G_1,G_1,\psi) \cdot \left ( p(1-\varrho)+\frac{\varrho}{2} \right ).
\end{eqnarray}
We will show in Section \ref{tjf} that, under the hypotheses of the first
 part of Theorem \ref{thaa}, $\Upsilon_R $ almost
always happens. Then, adding (\ref{tji}) to (\ref{tjj}) and using
(\ref{tjk}), we then obtain that almost always
$\Delta_R^{0+}(G_1,G_2,\psi) + \Delta_R^{0-}(G_1,G_2,\psi) < \frac{1}{2} \cdot
\Delta_R (G_1,G_1,\psi) $ for all $\psi \in \varPi_\phi\backslash\{e\}$.
Substituting this into (\ref{tjm}) yields
that almost always $\Delta_R (G_1,G_2,\psi)> \Delta_R (G_1,G_2,e)$
for all $\psi \in \varPi_\phi\backslash\{e\}$, and the first part
of Theorem \ref{thaa} will then be proven.

\subsection{Under Hypotheses of Theorem \ref{thaa}, Part i,
$\Upsilon_R $ Occurs Almost Always
\label{tjf} }

For any $k \in \{ 1,2,\ldots, m \}$, denote  $\varPi_\phi(k):=
\{ \psi \in \varPi_\phi: | \{ v \in V : \psi(v) \ne v \}  |=k \}$.
Just for now, let $k \in \{ 1,2,\ldots, m \}$ be chosen, and let $\psi
\in \varPi_\phi(k)$ be chosen. The random variable $\Delta_R(G_1,G_1,\psi)$ is a function
of the $\eta' :=  ks$ independent Bernoulli$(p)$
random variables 
$$ \{ \mathbbm{1}\{\{w,u\}\in E(G_1)\} \}_{ (w,u) \in W \times U:
\psi(w) \ne w},$$ and note that the hypotheses
of Theorem \ref{thg} are satisfied, hence for the choice of
$t=\frac{1}{20} \sqrt{\eta'p(1-p)}$ in Theorem \ref{thg} we obtain that
\begin{eqnarray} \label{tjx}
\mathbb{P} \left [ | \Delta_R(G_1,G_1,\psi)- \mathbb{E}\Delta_R(G_1,G_1,\psi) |>
\frac{1}{5} \eta'p(1-p) \right ]
\leq 2e^{- \eta' p (1-p)/400 }.
\end{eqnarray}
Also note that
\begin{eqnarray*}
\Delta_R(G_1,G_1,\psi)=\sum_{ \substack{(w,u) \in W \times U  \\ \text{s.t.} \psi(w) \ne w }}
 \mathbbm{1}\bigg\{ \mathbbm{1}\{\{ w,u \}\in E(G_1)\}  
\ne \mathbbm{1}\{ \{ \psi(w),u \}\in E(G_1) \}  \bigg\}
\end{eqnarray*}
is the sum of $\eta'$ Bernoulli$(2p(1-p))$ random variables hence
\begin{eqnarray} \label{tjy}
\mathbb{E} \Delta_R(G_1,G_1,\psi)=2 \eta'p(1-p) .
\end{eqnarray}
Thus, by (\ref{tjx}) and (\ref{tjy})
we obtain that
\begin{eqnarray} \mathbb{P} \left ( \frac{\Delta_R(G_1,G_1,\psi)}{ksp(1-p)} \not \in
[9/5,\  11/5] \right ) \leq 2 e^{\frac{-1}{400}ksp(1-p)}
\leq 2 e^{\frac{-\xi_2^2}{400}ks}  .  \label{tjp}
\end{eqnarray}

Conditioning on $G_1$, random variable $\Delta_R^{0+}(G_1,G_2,\psi)$ has a
$$\text{Binomial}\left (\Delta_R^+(G_1,G_1,\psi),(1-p)(1-\varrho) \right )$$ distribution, and
random variable $\Delta_R^{0-}(G_1,G_2,\psi)$ has a $$\text{Binomial}\left (
\Delta_R^-(G_1,G_1,\psi),p(1-\varrho) \right )$$distribution.
Conditioning also on the event that $\frac{\Delta_R(G_1,G_1,\psi)}{ksp(1-p)} \in
[9/5,\  11/5]$, applying Theorem \ref{thi} with the value $t=\frac{\varrho}{2}
\cdot \Delta_R^+(G_1,G_1,\psi)$, and using (\ref{tjk}), we have that
\begin{eqnarray} \label{tjn}
\mathbb{P} \left [  \Delta_R^{0+}(G_1,G_2,\psi) \geq
\Delta_R^+(G_1,G_1,\psi) \cdot \left (  (1-p)(1-\varrho)+ \frac{\varrho}{2} \right )
 \right ] \leq e^{\frac{-\xi_2^4}{20} \cdot ks}, \ \\ \label{tjo}
\mathbb{P} \left [  \Delta_R^{0-}(G_1,G_2,\psi) \geq
\Delta_R^-(G_1,G_1,\psi) \cdot \left (  p(1-\varrho)+ \frac{\varrho}{2} \right )
 \right ] \leq e^{\frac{-\xi_2^4}{20} \cdot ks}.
\end{eqnarray}

Finally, applying (\ref{tjp}), (\ref{tjn}) and (\ref{tjo}), the probability of $\Upsilon_R ^C$ can be bounded using subadditivity:
\begin{eqnarray*}
\mathbb{P}(\Upsilon_R ^C) & \leq & \sum_{k=1}^m \sum_{\psi \in \varPi_\phi(k)}
\left (  2 e^{\frac{-\xi_2^2}{400}ks} + e^{\frac{-\xi_2^4}{20} \cdot ks}  +  e^{\frac{-\xi_2^4}{20} \cdot ks} \right ) \\
 & \leq &  \sum_{k=1}^m m^k \left (  2 e^{\frac{-\xi_2^2}{400}ks} + 2e^{\frac{-\xi_2^4}{20} \cdot ks}  \right )\\
 & \leq & \sum_{k=1}^m \left (  2 e^{\frac{-\xi_2^2}{400}ks + k \log m} + 2e^{\frac{-\xi_2^4}{20} \cdot ks + k \log m}  \right )
 \leq m \cdot \frac{4}{m^3},
\end{eqnarray*}
the last inequality holding if
$s \geq c_5 \log m$  for sufficiently large, fixed
constant $c_5$. Because $\sum_{m=1}^\infty \frac{4}{n^2}<\infty$ we have by
the Borel-Cantelli Lemma that $\Upsilon_R $ almost always happens.
As mentioned in Section \ref{tjd}, this completes
the proof of the first part of Theorem \ref{thaa}.~$\blacksquare$
\begin{remark}
\emph{ We do not chase the exact constant $c_5$ here, focusing on the order of magnitude of $s$ instead.  Also, if we allow $p$ and $\rho$ to vary with $m$, then a minor alteration of the above proof (and a tighter Chernoff-Hoeffding bound) yields the same conclusion as in Theorem~\ref{thaa}.i  if for (an arbitrary but) fixed $0<\epsilon<2$ and $q_1:=(1-p)(1-\varrho)$ and $q_2:=p(1-\varrho)$
\begin{align*}
c_5&:=c_5(p,\varrho)\\
&>\max\left\{\frac{2}{H(q_1+\frac{\varrho}{2},q_1)\cdot p(1-p)(2-\epsilon)},\frac{2}{H(q_2+\frac{\varrho}{2},q_2)\cdot p(1-p)(2-\epsilon)},\frac{16}{\epsilon^2 p(1-p)}\right\}.
\end{align*}
Details are left to the reader.
}
\end{remark}

\subsection{Proof of the Theorem \ref{thaa}, Part ii\label{tka}}

We now prove the second part of Theorem \ref{thaa}.

Just for now, let $\psi \in \varPi(m)$ be chosen (i.e., none of the
nonseeds are fixed points for $\psi$),
and condition on $\Delta_R^+(G_1,G_1,\psi)=L$, where $ \frac{9}{10}smp(1-p) \leq L \leq
\frac{11}{10}smp(1-p)$. The random variables
$\Delta_R^{0+}(G_1,G_2,\psi)$ and $\Delta_R^{0-}(G_1,G_2,\psi)$
are independent, and have distributions Binomial$(L,q_1)$ and
Binomial$(L,q_2)$, respectively,
where $q_1:=(1-p)(1-\varrho)$ and $q_2:=p(1-\varrho)$.

Denoting $r_1:=q_1+\frac{\rho}{2}$ and $r_2:=q_2+\frac{\varrho}{2}$,
we have by Proposition \ref{P:bin} that
\begin{align*}
\mathbb{P} \Bigg ( \Delta_R^{0+}(G_1,G_2,\psi) &> L \cdot r_1 \ \mbox{ and } \
\Delta_R^{0-}(G_1,G_2,\psi) > L \cdot r_2 \Bigg )\\
 &\geq
\frac{q_1q_2}{200L}\sqrt{\frac{(1-r_1)(1-r_2)}{r_1r_2}}
e^{-L \cdot H(r_1,q_1) - L \cdot H(r_2,q_2)}
\end{align*}
Considering the bound on $H(x+y,y)$ described in Section \ref{spt},
we have that $H(r_1,q_1)$ and $H(r_2,q_2)$ are both bounded
above by a constant. With the fact that $1-r_1=r_2$ and
$1-r_2=r_1$, from the above we obtain that there is a positive
real number $c$ such that
\begin{align} \label{tkb}
\mathbb{P} \Bigg ( \Delta_R^{0+}(G_1,G_2,\psi) > L \cdot r_1 \ \mbox{ and }  \
\Delta_R^{0-}(G_1,G_2,\psi) > L \cdot r_2 \Bigg ) &\geq
\frac{c}{sm}\cdot e^{-\frac{sm}{c} }\notag\\
 &\geq \frac{c}{m \log m}\cdot e^{-\frac{sm}{c} }
\end{align}
 under the hypotheses of the second part of
Theorem \ref{thaa}.

Next,  $|\varPi_\phi(m)|$ is the number of derangements of an $m$ element set, and it
satisfies $\lim _{m \rightarrow \infty}
\frac{|\varPi_\phi (m)|}{m!}=\frac{1}{e}$,
thus with Stirling's formula we have that for $m$ large enough it will hold that $|\varPi(m)| \geq \left (  \frac{m}{e}   \right )^m$.
Thus, for $m$ large enough, by  (\ref{tjm}) and (\ref{tkb}),
\begin{align*}
\mathbb{E} |
\left \{ \psi \in \varPi_\phi : \Delta_R(G_1,G_2,\psi) < \Delta_R (G_1,G_2,e)
\right  \}| &= \sum_{\psi \in \varPi_\phi} \mathbb{P}
\Bigg (\Delta_R(G_1,G_2,\psi) < \Delta_R(G_1,G_2,e)  \Bigg )\\
& \geq  \sum_{\psi \in \varPi_\phi(m)} \mathbb{P}
\Bigg (\Delta_R (G_1,G_2,\psi) < \Delta_R(G_1,G_2,e)  \Bigg )\\
&\geq \left (  \frac{m}{e}   \right )^m  \frac{c}{m \log m}\cdot e^{-\frac{sm}{c} } \\
& =
\frac{c}{m \log m} \cdot e ^{-\frac{sm}{c} +m \log m - m},
\end{align*}
so that there exists a fixed real number $c_6>0$ such that
if $s \leq c_6 \log m$
then it follows that $\mathbb{E} |
\left \{ \psi \in \varPi_{\phi} : \Delta_R(G_1,G_2,\psi) < \Delta_R(G_1,G_2,e)
\right  \}| \rightarrow \infty $ as $m \rightarrow \infty$,
and Theorem~\ref{thaa} part ii~is~proven.~$\blacksquare$
\begin{remark}
\emph{
We could tighten the constant $c_6$ here, but choose instead to focus on the order of magnitude of $s$.  If we allow $p$ and $\varrho$ to be functions of $m$, then a simple alteration of the above proof yields the same results of Theorem~\ref{thaa}.ii  if
$$c_6:=c_6(p,\varrho)< \frac{1}{4\big[H(q_1+\frac{\varrho}{2},q_1)+H(q_2+\frac{\varrho}{2},q_2)\big]p(1-p)};$$
again details are left to the reader.
}
\end{remark}
\end{document}